\documentclass[moor]{informs2}              % for a regular run
\usepackage{natbib}
\NatBibNumeric
\def\bibfont{\small}%
\bibpunct[, ]{[}{]}{,}{n}{}{,}%

\TheoremsNumberedThrough     % Preferred (Theorem 1, Lemma 1, Theorem 2)
\EquationsNumberedThrough    % Default: (1), (2), ...
\usepackage{latexsym}
\usepackage{amsmath, amssymb, colonequals}
\usepackage{graphics, graphicx, url, color, setspace, lineno, enumerate, verbatim}
\usepackage{tikz}
\usepackage{pdflscape}
\usepackage{xspace}
\usepackage{bm}
\usepackage{comment}
\usepackage{float}
\usepackage{caption}
\usepackage{subcaption}
\usepackage{algorithmic}
\usepackage{algorithm}
\usepackage{pbox}
\usepackage[titletoc,toc,page]{appendix}
\usepackage{multirow}							%to add rows in the table that spans multiple rows
\usepackage[misc]{ifsym}						%for the envelope symbol
\usepackage{resizegather}
\usepackage{xr}

\usepackage{setspace}
\usepackage{authblk}

\def\tr{^{\intercal}}

\def\Re{{\mathbb R}}

\def\conv{\mathop{\rm conv}}

\def\Re{{\mathbb R}}

\def\ECR{EC\&R\xspace}

\newcommand{\vc}[1]{\bm{#1}}	%defines vector format
\newcommand{\mcl}[1]{\mathcal{#1}}	%defines font format
\newcommand{\mr}[1]{\mathrm{#1}}	%defines font format
\newcommand{\mt}[1]{\mathtt{#1}}	%defines font format

\begin{document}
	%%%%%%%%%%%%%%%%
	
	% Outcomment only when entries are known. Otherwise leave as is and 
	%   default values will be used.
	%\setcounter{page}{1}
	%\VOLUME{00}%
	%\NO{0}%
	%\MONTH{Xxxxx}% (month or a similar seasonal id)
	%\YEAR{0000}% e.g., 2005
	%\FIRSTPAGE{000}%
	%\LASTPAGE{000}%
	%\SHORTYEAR{00}% shortened year (two-digit)
	%\ISSUE{0000} %
	%\LONGFIRSTPAGE{0001} %
	%\DOI{10.1287/xxxx.0000.0000}%
	
	% Author's names for the running heads
	% Sample depending on the number of authors;
	% \RUNAUTHOR{Jones}
	% \RUNAUTHOR{Jones and Wilson}
	% \RUNAUTHOR{Jones, Miller, and Wilson}
	\RUNAUTHOR{Khademnia and Davarnia}
	% \RUNAUTHOR{Jones et al.} % for four or more authors
	% Enter authors following the given pattern:
	%\RUNAUTHOR{}
	
	% Title or shortened title suitable for running heads. Sample:
	% \RUNTITLE{Bundling Information Goods of Decreasing Value}
	% Enter the (shortened) title:
	\RUNTITLE{Convexification of Bilinear Terms over Network Polytopes}
	
	% Full title. Sample:
	% \TITLE{Bundling Information Goods of Decreasing Value}
	% Enter the full title:
	\TITLE{Convexification of Bilinear Terms over Network Polytopes}
	
	% Block of authors and their affiliations starts here:
	% NOTE: Authors with same affiliation, if the order of authors allows, 
	%   should be entered in ONE field, separated by a comma. 
	%   \EMAIL field can be repeated if more than one author
	\ARTICLEAUTHORS{%
		\AUTHOR{Erfan Khademnia}
		\AFF{Industrial and Manufacturing Systems Engineering, Iowa State University, Ames, IA, USA, \EMAIL{erfank@iastate.edu}, \URL{}}
		\AUTHOR{Danial Davarnia}
		\AFF{Industrial and Manufacturing Systems Engineering, Iowa State University, Ames, IA, USA, \EMAIL{davarnia@iastate.edu}, \URL{}}
		% Enter all authors
	} % end of the block
	\ABSTRACT{%
It is well-known that the McCormick relaxation for the bilinear constraint $z=xy$ gives the convex hull over the box domains for $x$ and $y$.
In network applications where the domain of bilinear variables is described by a network polytope, the McCormick relaxation, also referred to as linearization, fails to provide the convex hull and often leads to poor dual bounds.
We study the convex hull of the set containing bilinear constraints $z_{i,j}=x_iy_j$ where $x_i$ represents the arc-flow variable in a network polytope, and $y_j$ is in a simplex.
For the case where the simplex contains a single $y$ variable, we introduce a systematic procedure to obtain the convex hull of the above set in the original space of variables, and show that all facet-defining inequalities of the convex hull can be obtained explicitly through identifying a special tree structure in the underlying network. 
%This approach improves the classical McCormick bounds for network problems by taking advantage of the graphical structure of the network.
For the generalization where the simplex contains multiple $y$ variables, we design a constructive procedure to obtain an important class of facet-defining inequalities for the convex hull of the underlying bilinear set that is characterized by a special forest structure in the underlying network.
Computational experiments conducted on different applications show the effectiveness of the proposed methods in improving the dual bounds obtained from alternative techniques. 
	}%
	
	% Sample
	%\KEYWORDS{deterministic inventory theory; infinite linear programming duality; 
	%  existence of optimal policies; semi-Markov decision process; cyclic schedule}
	%\MSCCLASS{Primary: 90B05; secondary: 90C40, 90C90}
	%\ORMSCLASS{Primary: Inventory/production: deterministic multi-item;
	%  secondary: dynamic programming/optimal control: deterministic 
	%  semi-Markov; programming: infinite dimensional}
	%\HISTORY{Received November 20, 2003; revised March 8, 2004, and March 26, 2004.}
	
	% Fill in data. If unknown, outcomment the field
	%\MSCCLASS{}
	%\ORMSCLASS{Primary: ; secondary: }
	%\HISTORY{}
	
	\maketitle
	%%%%%%%%%%%%%%%%%%%%%%%%%%%%%%%%%%%%%%%%%%%%%%%%%%%%%%%%%%%%%%%%%%%%%%
	
	% Samples of sectioning (and labeling) in MOOR.
	% NOTE: (1) all section levels end with a period,
	%       (2) capitalization is as shown (sentence style, not title style).
	%
	%\section{Introduction.}\label{intro} %%1.
	%\subsection{Duality and the classical EOQ problem.}\label{class-EOQ} %% 1.1.
	%\subsection{Outline.}\label{outline1} %% 1.2.
	%\subsubsection{Cyclic schedules for the general deterministic SMDP.}
	%  \label{cyclic-schedules} %% 1.2.1
	%\section{Problem description.}\label{problemdescription} %% 2.
	
	% Text of your paper here

%********************************************************************************************************************************************************
%********************************************************************************************************************************************************
%********************************************************************************************************************************************************

% Insert the name of "your journal" with
%\journalname{}

%%%%%%%%%%%%%%%%%%%%%%%
%%%%%%%%%%%%%%%%%%%%%%%%%%%%%%%%%%%%%%%%%%%%%%
%%%%%%%%%%%%%%%%%%%%%%%%%%%%%%%%%%%%%%%%%%%%%%%%%%%%%%%%%%%%%%%%%%%%%
%%%%%%%%%%%%%%%%%%%%%%%%%%%%%%%%%%%%%%%%%%%%%%%%%%%%%%%%%%%%%%%%%%%%%%%%%%%%%%%%%%%%%%%%%%%%
\section{Introduction} \label{sec:introduction}
%%%%%%%%%%%%%%%%%%%%%%%%%%%%%%%%%%%%%%%%%%%%%%%%%%%%%%%%%%%%%%%%%%%%%%%%%%%%%%%%%%%%%%%%%%%%
%%%%%%%%%%%%%%%%%%%%%%%%%%%%%%%%%%%%%%%%%%%%%%%%%%%%%%%%%%%%%%%%%%%%%
%%%%%%%%%%%%%%%%%%%%%%%%%%%%%%%%%%%%%%%%%%%%%%
%%%%%%%%%%%%%%%%%%%%%%%

Bilinear constraints in conjunction with network models appear in various mixed-integer and nonlinear programming (MINLP) applications, including the decision diagrams \cite{davarnia:va:2020,davarnia:2021,salemi:davarnia:2022,salemi:da:2023}, the fixed-charge network flow problems \cite{rebennack:na:pa:2009}, network models with complementarity constraints, such as transportation problems with conflicts \cite{ficker:sp:wo:2021} and the red-blue transportation problems \cite{vancroonenburg:de:go:sp:2014}. 
A common occurrence of bilinear terms in network models pertains to bilevel network problems after being reformulated as a single-level program through either using a dual formulation of the inner problem or incorporating optimality conditions inside the outer problem \cite{benayed:bl:1990,chiou:2005}.
These reformulation approaches are widely used in the network interdiction problems, where newly added bilinear terms are relaxed using a \textit{linearization} technique based on the McCormick bounds \cite{mccormick:1976} over a box domain; see \cite{smith:li:2008} for an exposition.
While these relaxations provide the convex hull of the bilinear constraint over a box domain \cite{alkhayyal:fa:1983}, they often lead to weak relaxations when the variables domain becomes more complicated as in general polyhedra \cite{gupte:ah:ch:de:2013,davarnia:2016}.
It has been shown \cite{boland:de:ka:mo:ri:2017,luedtke:na:li:2012} that when the number of bilinear terms in the underlying function increases, while the variables are still in a box domain, the McCormick bounds can be very poor compared to the ones obtained from the convex and concave envelopes of the function. 
It is shown in \cite{bonami:2018} that while McCormick relaxation performs poorly for box-constrained quadratic problems, a stronger relaxation can be achieved by adding Chv\'{a}tal-Gomory valid inequalities.

\medskip
There are various studies in the literature that develop convexification methods for multilinear functions as a generalization of bilinear forms, but the side constraints for the involved variables are often limited to variable bounds; see \cite{delpia:kh:2018,delpia:kh:2021,khajavirad:2023} for an exposition to different classes of valid inequalities for multilinear sets.
%For instance, \cite{delpia:kh:2021} introduces a new class of valid inequalities, called running intersection inequalities, for the multilinear polytope described by a set of multilinear equations.
In \cite{delpia:di:2021,delpia:wa:2022}, the authors study the Chv\'{a}tal rank of different classes of cutting planes, including those associated with $\beta$-cylces of the underlying hypergraph, in the multilinear polytope that arises in binary polynomial problems, while \cite{buchheim:cr:ro:2019} provides a characterization of the multilinear polytope in terms of the acyclicity of its corresponding hypergraph, which guarantees integrality of optimal solutions obtained by the standard linearization procedure.
The authors in \cite{gupte:ka:ri:wa:2020} derive extended formulations for the convex hull of the graph of a bilinear function on the $n$-dimensional unit cube through identifying the facets of the Boolean Quadratic Polytope.
In \cite{fampa:lee:2021}, an efficient method to convexify bilinear functions through McCormick relaxations is proposed which takes advantage of the structural convexity in a symmetric quadratic form.
Further, \cite{muller:se:gl:2020} introduces a new class of valid inequalities, tailored for the spatial branch-and-bound framework, based on projection techniques that improve the quality of McCormick inequalities.
Other works in the literature consider a polyhedral, often triangular, subdivision of the domain to derive strong valid inequalities for a bilinear set; see \cite{tawarmalani:ri:xi:2012,sherali:al:1992,locatelli:sch:2014} for examples of such approaches.
Further, \cite{davarnia:ri:ta:2017} proposes an aggregation procedure, referred to as \textit{extended cancel-and-relax} (EC\&R), to simultaneously convexify the graph of bilinear functions over a general polytope structure. 
In this paper, we make use of the \ECR procedure to derive convexification results for a bilinear set where the side constraints on variables are described by a network flow model as defined next.

\medskip
For $N := \{1,\dotsc,n\}$, $M := \{1,\dotsc,m\}$, $K := \{1,\dotsc,\kappa\}$, and $T := \{1, \dotsc, \tau\}$, we consider
\[
\mcl{S} = \left\{ (\vc{x};\vc{y};\vc{z}) \in \Xi \times \Delta_m \times \Re^{\kappa} \, \middle|\,
\vc{y}\tr A^k \vc{x} = z_{k}, \, \, \, \forall k \in K
\right\},
\]
where $\Xi = \left\{ \vc{x} \in \Re^n \, \middle| \, E\vc{x} \geq \vc{f}, \, \vc{0} \leq \vc{x} \leq \vc{u} \right\}$ is a \textit{primal} network polytope described by the flow-balance and arc capacity constraints, and $\Delta_m = \left\{ \vc{y} \in \Re_+^{m} \, \middle| \, \vc{1}\tr \vc{y} \leq 1 \right\}$ is a simplex.
When variables $\vc{y}$ are binary, $\Delta_m$ represents a \textit{special ordered set of type I} (SOS1); see \cite{beale:fo:1976} for an exposition.
Such simplex structures appear in various applications and can be obtained by reformulating the underlying polytopes through extreme point decomposition; see \cite{davarnia:ri:ta:2017} for a detailed account. 
In the above definition, $E \in \Re^{\tau \times n}$, $\vc{f} \in \Re^{\tau}$, $\vc{u} \in \Re^n$, and $A^k \in \Re^{m\times n}$ is a matrix with all elements equal to zero except one that is equal to one, i.e., if $A^k_{ji} = 1$ for some $(i,j) \in N \times M$, the bilinear constraint with index $k$ represents $y_jx_i = z_k$.
%We refer to the $j^{\textrm{th}}$ row (resp. $j^{\textrm{th}}$ column) of $A^k$ by $A^k_{j.}$ (resp. by $A^k_{.j}$).

\medskip
The contributions of this paper are as follows.
We propose a systematic procedure to convexify $\mcl{S}$ and derive explicit inequalities in its description.
The resulting cutting planes are directly obtained in the original space of variables.
We show that facet-defining inequalities in the convex hull description can be explicitly derived by identifying special tree and forest structures in the underlying network, leading to an interpretable and efficient cut-generating oracle.   
The inequalities obtained from our proposed algorithms can be added to the typical McCormick relaxations to strengthen the formulation and improve the bounds. 

\medskip
Since we use the specialized aggregation technique introduced in \cite{davarnia:ri:ta:2017} in our derivation process, we next explain the novelties of this work compared to \cite{davarnia:ri:ta:2017}. 
In particular, \cite{davarnia:ri:ta:2017} applies the proposed aggregation technique to a special case of $\mcl{S}$ where $\Xi$ is a \textit{dual} network polytope described by the dual of the network flow formulation, inspired by the network interdiction applications, with the goal of obtaining explicit forms for the facet-defining inequalities that would be practical to generate. 
These \textit{explicit} convexification results are obtained only for the case with $m=1$ and $k=1$. 
In this paper, we extend and generalize these results from two directions.
(i) We obtain the convexification results for set $\mcl{S}$ where $\Xi$ is a \textit{primal} network polytope, which appears in a broader application domain beyond the network interdiction problem, such as those discussed above. 
Due to the structural difference between the primal and dual network polytopes, the framework developed in \cite{davarnia:ri:ta:2017} is not applicable to the network problems studied in this paper, advocating the need for a novel framework specifically designed for this new set. 
Consequently, the form of facet-defining inequalities for the network polytope obtained in this paper is fundamentally different from those derived in \cite{davarnia:ri:ta:2017}.
(ii) We obtain explicit convexification results for cases with general $m$ and $k$. 
This generalization makes the problem structure, and hence the analysis, much more complicated due to losing the totally unimodularity property that exists when $m = k = 1$; see Section~\ref{sec:primal} for a detailed discussion. From a practical point of view, as shown in Section~\ref{sec:computation}, the cutting planes obtained through our newly designed algorithms for the general cases with $m > 1$ and $k > 1$ can significantly improve the bounds obtained by adding cutting planes for the cases with $m = k = 1$ as well as the bounds obtained by alternative methods such as reformulation-linearization technique (RLT), highlighting the practical advantage and importance of the extended convexification results developed in this paper.

%The presented methods consider a general network structure, which complement and generalize the results of \cite{davarnia:ri:ta:2017} that are obtained for network interdiction problems.
%In particular, \cite{davarnia:ri:ta:2017} presents the convexification results for a special case of $\mcl{S}$ where $m = 1$, $\kappa = 1$ and $\Xi$ is a dual network polytope.
%In this work, we extend these results by considering the cases where $M$ and $K$ can have multiple elements and $\Xi$ is a primal network polytope.

\medskip
The remainder of the paper is organized as follows.
We give a brief background on the \ECR procedure as a basis of our analysis in Section~\ref{sec:ECR}.
In Section~\ref{sec:primal}, we obtain convexification results for bilinear terms defined over a network polytope.
Preliminary computational experiments are presented in Section~\ref{sec:computation} to show the effectiveness of the developed cut-generating frameworks.
We give concluding remarks in Section~\ref{sec:conclusion}.

\medskip
\textbf{Notation.} 
Bold letters represent vectors.
We refer to the $j^{\textrm{th}}$ row (resp. $j^{\textrm{th}}$ column) of a matrix $A \in \Re^{m \times n}$ by $A_{j.}$ (resp. by $A_{.j}$).
For a given set $S \subseteq \Re^n$, we denote by $\conv(S)$ its convex hull.
We use symbol $\pm$ to show both cases with $+$ and $-$. 
For example, when we use $l^{\pm}$ in an expression, it means that expression holds for both cases $l^+$ and $l^-$.

%%%%%%%%%%%%%%%%%%%%%%%
%%%%%%%%%%%%%%%%%%%%%%%%%%%%%%%%%%%%%%%%%%%%%%
%%%%%%%%%%%%%%%%%%%%%%%%%%%%%%%%%%%%%%%%%%%%%%%%%%%%%%%%%%%%%%%%%%%%%
%%%%%%%%%%%%%%%%%%%%%%%%%%%%%%%%%%%%%%%%%%%%%%%%%%%%%%%%%%%%%%%%%%%%%%%%%%%%%%%%%%%%%%%%%%%%
\section{Extended Cancel-and-Relax} \label{sec:ECR}
%%%%%%%%%%%%%%%%%%%%%%%%%%%%%%%%%%%%%%%%%%%%%%%%%%%%%%%%%%%%%%%%%%%%%%%%%%%%%%%%%%%%%%%%%%%%
%%%%%%%%%%%%%%%%%%%%%%%%%%%%%%%%%%%%%%%%%%%%%%%%%%%%%%%%%%%%%%%%%%%%%
%%%%%%%%%%%%%%%%%%%%%%%%%%%%%%%%%%%%%%%%%%%%%%
%%%%%%%%%%%%%%%%%%%%%%%

In this section, we present the \ECR procedure adapted for $\mcl{S}$.
The following theorem is at the core of this procedure, as shown in Theorem 2.7 of \cite{davarnia:ri:ta:2017}.

\begin{theorem}
	A convex hull description of $\mcl{S}$ can be obtained by the linear constraints in $\Xi$ and $\Delta_m$ together will all class-$l^{\pm}$ \ECR inequalities for all $l \in K$. 
%	\Halmos
\end{theorem}	

A class-$l^{\pm}$ \ECR inequality is obtained through a weighted aggregation of the constraints in the description of $\mcl{S}$.
The procedure to generate this inequality is as follows.

\begin{enumerate}
	\item
We select $l \in K$ to be the index of a bilinear constraint used in the aggregation, which we refer to as the \textit{base} equality.
	We also select a sign indicator $+$ or $-$ to indicate whether a weight $1$ or $-1$ is used for the base equality during aggregation.
	\item
	Defining $T$ as the index set of the non-bound constraints in $\Xi$, we select $\mcl{I}_1$, $\cdots$, $\mcl{I}_m$ and $\bar{\mcl{I}}$ as subsets of $T$ whose intersection is empty. 
	Then, for each $j \in M$ and for each $t \in \mcl{I}_j$ (resp. $t \in \bar{\mcl{I}}$), we multiply the constraint $E_{t.}\vc{x} \geq f_t$ by $\gamma^j_ty_j$ where $\gamma^j_t \geq 0$ (resp. by $\theta_t(1-\sum_{i \in M}y_i)$ where $\theta_t \geq 0$).
	\item
	We select $\mcl{J}$ and $\bar{\mcl{J}}$ as disjoint subsets of $N$. 
	Then, for each index $i \in \mcl{J}$, we multiply $x_i \geq 0$ by $\lambda_i(1-\sum_{j \in M}y_j)$ where $\lambda_i \geq 0$, and for each $i \in \bar{\mcl{J}}$, we multiply $u_i - x_i \geq 0$ by $\mu_i(1-\sum_{j \in M}y_j)$ where $\mu_i \geq 0$.
\end{enumerate}
The above sets are compactly represented as $\big[\mcl{I}_1,\dotsc,\mcl{I}_m,\bar{\mcl{I}} \big| \mcl{J},\bar{\mcl{J}}\big]$, which is called an \textit{\ECR assignment}.
Each \ECR assignment is identified by its class-$l^{\pm}$ where $l$ is the index of the base equality and $\pm$ is its sign indicator.
We next aggregate all aforementioned weighted constraints.
During the aggregation, we require that weights $\gamma$, $\theta$, $\lambda$ and $\mu$ be chosen in such a way that:
\begin{itemize}
\item[\textit{(C1)}] 
	at least $\sum_{j \in M}|\mcl{I}_j|+|\bar{\mcl{I}}|+|\mcl{J}|+|\bar{\mcl{J}}|$ bilinear terms are canceled (i.e., their coefficient becomes zero), and
	\item[\textit{(C2)}]
	if $\cup_{j \in M}\mcl{I}_j \cup \bar{\mcl{I}} \cup \mcl{J} \cup \bar{\mcl{J}} \neq \emptyset$, for each constraint used in the aggregation (including the base equality), at least one bilinear term among all those created after multiplying that constraint with its corresponding weight is canceled.

\end{itemize}
The desired \ECR inequality is then obtained by relaxing (i.e., replacing) the remaining bilinear terms of the form $x_iy_j$ in the aggregated inequality using either $u_iy_j - x_iy_j \geq 0$ or $-x_iy_j + z_k \geq 0$ for some $k \in K$ such that $A^k_{ji} = 1$. Similarly, the bilinear terms of the form $-x_iy_j$ is relaxed using either $x_iy_j \geq 0$ or $x_iy_j - z_k \geq 0$ for some $k \in K$ such that $A^k_{ji} = 1$.
The resulting linear inequality is referred to as a class-$l^{\pm}$ \ECR inequality.

\medskip
A summary of the derivation steps for the \ECR procedure is given in the appendix; see \cite{davarnia:ri:ta:2017} for a detailed exposition.
Here, we only present an end result that will be used in the subsequent sections.
%general form of a constraint set whose feasible solutions correspond to the \ECR inequalities; see equation~\eqref{eq:Cl} below.}

\medskip
For each $l \in K$, define $K_l = K \setminus \{l\}$, and let $\vc{\pi}^l = \left(\vc{\beta}^+;\vc{\beta}^-;\{\vc{\gamma}^j\}_{j \in M};\vc{\theta};\{\vc{\eta}^j\}_{j \in M};\{\vc{\rho}^j\}_{j \in M};\vc{\lambda};\vc{\mu}\right)$ be the vector of aggregation weights used in the \ECR procedure.
In particular, the index $l$ is determined in Step 1 of the \ECR procedure by picking a base equality $l$ with either $+1$ or $-1$ weights.
The components $\vc{\gamma}^j$ (resp. $\vc{\theta}$) capture the weights for $y_j$ (resp. $1 - \sum_{j \in M} y_j$) when multiplied with the non-bound constraints in $\Xi$ as demonstrated in Step 2 of the \ECR procedure.
Similarly, $\vc{\lambda}$ and $\vc{\mu}$ record the weights for $1 - \sum_{j \in M} y_j$ when multiplied with the bound constraints in $\Xi$ in Step 3 of the \ECR procedure.
Finally, the relaxation step in the \ECR procedure determines the components $\vc{\beta}^{\pm}$ as the weights of the bilinear constraints in $K_l$, as well as the components $\vc{\eta}^j$ and $\vc{\rho}^j$ as the weights for $y_j$ when multiplied with the bound constraints in $\Xi$ to cancel the remaining bilinear terms.
Then, it is shown in Theorem 2.7 of \cite{davarnia:ri:ta:2017} that criteria (C1) and (C2) of the \ECR procedure provide necessary conditions for the weight vector $\vc{\pi}^l$ to be an extreme point of
\begin{multline} \label{eq:Cl}
	\mathcal{C}^{l} =
	\left\{
	\vc{\pi}^l \in \Re^{2(\kappa-1) + (m+1)(\tau + 2n)}_{+} \, \middle|
	\sum_{k \in K_l}A^k_{ji}\left(\beta^+_k - \beta^-_k\right) + \sum_{t \in T} E_{ti} \left(\gamma^j_t - \theta_t\right) \right. \\
	\left. + \eta^j_i - \rho^j_i - \lambda_i + \mu_i = \pm A^l_{ji}, \, \forall (i,j) \in N \times M \right\}.
\end{multline}

This set plays a central role in our analysis because the collection of the \ECR inequalities corresponding to the extreme points of $\mathcal{C}^{l}$, for all $l \in K$, contains all \textit{non-trivial} facet-defining inequalities in $\conv(\mcl{S})$, where a non-trivial inequality is one that cannot be implied by the linear constraints in the description of $\Xi$ and $\Delta_m$.

\section{Network Polytopes} \label{sec:primal}
%%%%%%%%%%%%%%%%%%%%%%%%%%%%%%%%%%%%%%%%%%%%%%%%%%%%%%%%%%%%%%%%%%%%%%%%%%%%%%%%%%%%%%%%%%%%
%%%%%%%%%%%%%%%%%%%%%%%%%%%%%%%%%%%%%%%%%%%%%%%%%%%%%%%%%%%%%%%%%%%%%
%%%%%%%%%%%%%%%%%%%%%%%%%%%%%%%%%%%%%%%%%%%%%%
%%%%%%%%%%%%%%%%%%%%%%%

In this section, we study the set $\mcl{S}$ where $\Xi$ represents a network polytope.
In particular, the constraint set $E\vc{x} \geq \vc{f}$ is composed of the flow-balance constraints after being separated into two inequalities of opposite signs, i.e., $E$ is an augmented node-arc incidence matrix where each row is duplicated with a negative sign, and $\vc{f}$ represents the extended supply/demand vector.
In this description, $\vc{u}$ denotes the arc-capacity vector.
\subsection{The Case with $m = 1$.} \label{subsec:primal-single}
%%%%%%%%%%%%%%%%%%%%%%%%%%%%%%%%%%%%%%%%%%%%%
%%%%%%%%%%%%%%%%%%%%%%%%%%%%%%%%%%%%%%%%%
%%%%%%%%%%%%%%%%%%%%%%%%%%%%%%%%

In this section, we consider the case where $m = 1$, whose corresponding bilinear set is denoted by $\mcl{S}^1$.
This special case is important as it provides a natural setting to improve the McCormick results as the most common relaxation for bilinear programs, which considers $m = 1$, even when there are multiple $\vc{y}$ variables that interact with each other through side constraints. 
In fact, the McCormick relaxation of bilinear terms with $m = 1$ is a building block of the factorable decomposition technique at the core of MINLP solvers \cite{mccormick:1976}. 
Therefore, the convexification results for the case where each $\vc{y}$ variable is treated separately (i.e., $m = 1$) is critical in determining the quality of the relaxation bound obtained for the entire model. 
For this reason, improving the McCormick relaxation, even for the simplest case with $m = 1$, could lead to substantial bound improvement as corroborated by the computational results in the Section \ref{sec:computation}.

\medskip
In this special case, we can simplify notation by matching the indices of $\vc{z}$ and $\vc{x}$ variables such that $y_1x_k = z_k$ for all $k \in K = N$.
We next show that, to generate an \ECR inequality for $\mcl{S}^1$, it is sufficient to use aggregation weight $1$ for all constraints used in the aggregation.

\begin{proposition} \label{prop:primal-weight}
	Let $\bar{\vc{\pi}}^l$ be an extreme point of the projection cone $\mcl{C}^l$, for some $l \in K$, corresponding to a non-trivial facet-defining inequality of $\conv(\mcl{S}^1)$. Then, it can be scaled in such a way that all of its components are 0 or 1. 
\end{proposition}	

\proof{Proof.}
	When $m=1$, we can write $\mcl{C}^l$ in \eqref{eq:Cl} as
	\begin{multline*}
		\mathcal{C}^{l} =
		\left\{
		\vc{\pi}^l \in \Re^{2(\kappa-1) + 2(\tau + 2n)}_{+} \, \middle|
		\sum_{t \in T} E_{ti} \left(\gamma^1_t - \theta_t\right) + \mu_i - \lambda_i \right. \\
		\left. + \eta^1_i - \rho^1_i  + \sum_{k \in K_l}A^k_{1i}\left(\beta^+_k - \beta^-_k\right) = \pm A^l_{1i}, \, \forall i \in N \right\}.
	\end{multline*}
	We can rearrange the columns of the coefficient matrix of the system defining $\mathcal{C}^{l}$ to obtain
	\begin{equation}
		\Bigl[
		\begin{array}{c|c|c|c|c|c|c|c}
			E^{\tr} \, & \, -E^{\tr} \, & \, I \, & \, -I \, & \, I \, & \, -I & \, \bar{I} \, & \, -\bar{I}
		\end{array}
		\Bigr]. \label{eq:ECR-matrix}
	\end{equation}
%	In the above matrix, the rows correspond to bilinear terms $y_1x_i$ (i.e., $w^1_i$ in the disjunctive programming formulation \eqref{eq_ef}) for $i \in N$. 
	In the above matrix, the rows correspond to bilinear terms $y_1x_i$ for $i \in N$. 
	The first and second column blocks correspond to the weights of the non-bound constraints in $\Xi$ multiplied by $y_{1}$ and $1-y_{1}$, which are denoted by $\gamma^1_t$ and $\theta_t$, respectively, for all $t \in T$.
	The third and fourth column blocks correspond to the weights of the lower and upper bound constraints on variables in $\Xi$ multiplied by $1-y_{1}$, which are captured by $\mu_i$ and $\lambda_i$, respectively, for all $i \in N$.
	Similarly, the fifth and sixth column blocks correspond to the weights of the lower and upper bound constraints on variables in $\Xi$ multiplied by $y_{1}$, which are recorded by $\eta^1_i$ and $\rho^1_i$, respectively, for all $i \in N$.
	In these columns, $I$ represents the identity matrix of appropriate size.
	Lastly, the seventh and eighth column blocks correspond to the weights of the bilinear constraints in $\mcl{S}^1$, which are represented by $\beta^+_k$ and $\beta^-_k$, respectively, for all $k \in K_l$.
	In particular, the element at column $k \in K_l$ and row $i \in N$ of $\bar{I}$ is equal to $1$ if $i=k$, and it is equal to zero otherwise.
	Based on these column values, it can be easily verified that \eqref{eq:ECR-matrix} is totally unimodular (TU).
	In $\mathcal{C}^{l}$, the right-hand-side vector is $\pm\vc{e}^{l} \in \Re^{m+n}$, where $\vc{e}^{l}$ is the unit vector whose components are all zero except for that corresponding to row $l$ representing $y_1x_l$, which is equal to $1$.
	Because $\bar{\vc{\pi}}^l$ is an extreme point of $\mathcal{C}^{l}$, it is associated with a basic feasible solution for its system of equations.
	Let $B$ be the corresponding basis for \eqref{eq:ECR-matrix}.
	It follows from Cramer's rule that all elements of $B^{-1}$ belong to $\{0,-1,1\}$ since \eqref{eq:ECR-matrix} is TU.
	Therefore, the components of $\pm B^{-1}\vc{e}^{l}$ belong to $\{0,-1,1\}$.
	We conclude that the components of basic feasible solutions to $\mathcal{C}^{l}$ are equal to $0$ or $1$ due to non-negativity of all variables in its description.
	\Halmos
\endproof	

\medskip
\begin{remark}	\label{rem:1-simplex}
	When $m = 1$, multiplying the bound constraints with $1-y_1$ in Step 3 of the \ECR procedure produces two of the standard McCormick bounds. 
	As a result, we can skip Step 3 in the \ECR procedure and merge it into the relaxation step, in which the other two McCormick bounds and the bilinear constraints are used for relaxing the remaining bilinear terms.
	In specific, any remaining bilinear term in the aggregated inequality can be \textit{relaxed} into either of the two McCormick lower bounds or the two McCormick upper bounds or the $\pm z$ variable corresponding to that term depending on its sign.
	In this case, the characterization of \ECR assignment can be reduced further to $\big[\mcl{I}_1,\bar{\mcl{I}}\big]$.  
\end{remark}

\medskip
\begin{remark}	\label{rem:1-separation}
	As described in Remark~\ref{rem:1-simplex}, each remaining bilinear term in the aggregated inequality of the \ECR procedure can be relaxed into three different linear terms. 
	While this can lead to an exponential growth in the number of resulting linear \ECR inequalities for each \ECR assignment, we can use an efficient separation procedure to find the most violated inequality among the resulting \ECR inequalities as follows.
	Assume that we aim to separate a given solution $(\bar{\vc{x}}; \bar{y}_1; \bar{\vc{z}})$ from $\conv(\mcl{S}^1)$ through the \ECR inequalities obtained from the aggregated inequality $g(\bar{\vc{x}}; y_1; \bar{\vc{z}}) \geq 0$ associated with the \ECR assignment $\big[\mcl{I}_1,\bar{\mcl{I}}\big]$.
	For each bilinear term $y_1x_i$, we choose the relaxation option that provides the minimum value among $u_i\bar{y}_1$ obtained from using $y_1(u_i - x_i) \geq 0$, $\bar{x}_i$ obtained from using $(1-y_1)x_i \geq 0$, and $\bar{z}_i$ obtained from using $-y_1x_i + z_i \geq 0$.
	Similarly, for each bilinear term $-y_1x_i$, we choose the relaxation option that provides the minimum value among $0$ obtained from using $y_1x_i \geq 0$, $u_i - \bar{x}_i - u_i\bar{y}_1$ obtained from using $(1-y_1)(u_i - x_i) \geq 0$, and $-\bar{z}_i$ obtained from using $y_1x_i - z_i \geq 0$.
	This approach provides the most violated \ECR inequality in the time linear in the number of remaining bilinear terms in the aggregated inequality.
\end{remark}

\medskip
Considering the relation between the extreme points of the projection cone $\mcl{C}^l$ for $l \in K$ and the aggregation weights in the \ECR procedure, Proposition~\ref{prop:primal-weight} and Remark~\ref{rem:1-simplex} imply that generating class-$l^{\pm}$ \ECR inequalities reduces to identifying the assignment $\big[\mcl{I}_1,\bar{\mcl{I}}\big]$ as the aggregation weights are readily determined.
In particular, the constraints in $\mcl{I}_1$ are multiplied with $y_1$, and those in $\bar{\mcl{I}}$ are multiplied with $1-y_1$.
We next show that, for set $\mcl{S}^1$, identifying all the \ECR assignments that satisfy the \ECR conditions (C1) and (C2) can be achieved by considering a special graphical structure in the underlying network.

\medskip
Given a network $\mr{G} = (\mr{V},\mr{A})$ with a node set $\mr{V}$ and arc set $\mr{A}$, assume that the index $k$ of variables $z_k$ in the description of $\mcl{S}^1$ refers to the arc whose flow variable $x_k$ appears in that bilinear constraint, i.e., $y_1x_k = z_k$ for $k \in \mr{A} = N = K$.
We define $t(k)$ and $h(k)$ to be the tail and head nodes of arc $k \in \mr{A}$, respectively.
Further, for any node $i \in \mr{V}$, we define $\delta^+(i)$ and $\delta^-(i)$ to be the set of outgoing and incoming arcs at that node, respectively.
We refer to the flow-balance inequality $\sum_{k \in \delta^+(i)} x_k - \sum_{k \in \delta^-(i)} x_k \geq f_i$ (resp. $-\sum_{k \in \delta^+(i)} x_k + \sum_{k \in \delta^-(i)} x_k \geq -f_i$) corresponding to node $i$ as the \textit{positive} (resp. \textit{negative}) flow-balance inequality, and refer to its index in the description of $\Xi$ by $i^+$ (resp. $i^-$) to be recorded in the \ECR assignment.     
For example, an \ECR assignment $\big[\{i^+\},\{j^-\}\big]$ implies that, in the aggregation, 
%the bilinear constraint corresponding to the arc with index $k \in \mr{A}$ is multiplied with $1$, the bilinear constraint corresponding to the arc with index $l \in \mr{A}$ is multiplied with $-1$, 
the positive flow-balance inequality corresponding to the node $i \in \mr{V}$ is multiplied with $y_1$, and the negative flow-balance inequality corresponding to the node $j \in \mr{V}$ is multiplied with $1-y_1$.
In the sequel, we denote the undirected variant of a subnetwork $\mr{P}$ of $\mr{G}$ by $\bar{\mr{P}}$, and conversely, we denote the the directed variant of an undirected subnetwork $\bar{\mr{P}}$ of $\bar{\mr{G}}$ by $\mr{P}$ according to the arc directions in $\mr{G}$. 

\begin{proposition} \label{prop:primal-tree}
	Consider set $\mcl{S}^1$ with $\Xi$ that represents the network polytope corresponding to the network $\mr{G} = (\mr{V},\mr{A})$.
	Let $\big[\mcl{I}_1,\bar{\mcl{I}}\big]$ be an \ECR assignment for class-$l^{\pm}$, for some $l \in \mr{A}$, that leads to a non-trivial facet-defining inequality of $\conv(\mcl{S}^1)$. 
	Define $\widetilde{\mcl{I}} = \{i \in \mr{V} | i^{\pm} \in \mcl{I}_1 \cup \bar{\mcl{I}} \}$ to be the subset of nodes whose flow-balance inequalities are used in the aggregation.
	Then, there exists a tree $\bar{\mr{T}}$ of $\bar{\mr{G}}$ composed of the nodes in $\widetilde{\mcl{I}}$ such that arc $l$ is incident to exactly one node of $\bar{\mr{T}}$.
\end{proposition}	

\proof{Proof.}
	First, we observe that for each node $i \in \mr{V}$, both of its positive and negative flow-balance inequalities cannot be selected for the aggregation, since otherwise, the columns representing the positive and negative inequalities in the basis of the coefficient matrix \eqref{eq:ECR-matrix} associated with the extreme point of $\mcl{C}^l$ would be linearly dependent, which would be a contradiction to the fact that the selected \ECR assignment leads to a facet-defining inequality of $\conv(\mcl{S}^1)$; see the proof of Proposition~\ref{prop:primal-weight} for details.
	As a result, considering that $\mcl{I}_1 \cap \bar{\mcl{I}} = \emptyset$ by the \ECR requirement, at most one of the following possibilities can occur in the \ECR assignment: $i^+ \in \mcl{I}_1$, $i^+ \in \bar{\mcl{I}}$, $i^- \in \mcl{I}_1$, and $i^- \in \bar{\mcl{I}}$.
	Therefore, each node in $\widetilde{\mcl{I}}$ corresponds to a unique flow-balance constraint in the \ECR assignment.
	Next, we show that arc $l$ is incident to exactly one node of $\widetilde{\mcl{I}}$.
	It follows from condition (C2) of the \ECR procedure that the bilinear term $y_1x_l$ for arc $l$ in the base equality must be canceled during the aggregation. 
	The constraints of $\Xi$ that can produce the bilinear term $y_1x_l$ during the aggregation are the flow-balance constraint corresponding to the tail node $t(l)$ of arc $l$, and the flow-balance constraint corresponding to the head node $h(l)$ of arc $l$.
	Since the aggregation weight for all the constraints in the \ECR assignment are $1$ according to Proposition~\ref{prop:primal-weight}, and considering that each flow-balance constraint can appear once in the aggregation as noted above, the only possibility to cancel the term $y_1x_l$ is to pick exactly one of the above constraints in the \ECR assignment.
	As a result, exactly one of the head and the tail nodes of arc $l$ must be in $\widetilde{\mcl{I}}$.
	Next, we show that there exists a tree $\bar{\mr{T}}$ of $\bar{\mr{G}}$ whose node set is $\widetilde{\mcl{I}}$.
	Assume by contradiction that there is no such tree composed of the nodes in $\widetilde{\mcl{I}}$.
	Therefore, $\widetilde{\mcl{I}}$ can be partitioned into two subsets $\widetilde{\mcl{I}}_1$ and $\widetilde{\mcl{I}}_2$, where the nodes in $\widetilde{\mcl{I}}_1$ are not adjacent to any nodes in $\widetilde{\mcl{I}}_2$.
	It is clear that arc $l$ cannot be incident to the nodes in both $\widetilde{\mcl{I}}_1$ and $\widetilde{\mcl{I}}_2$, since otherwise $\widetilde{\mcl{I}}_1$ and $\widetilde{\mcl{I}}_2$ would have adjacent nodes.
	Assume without the loss of generality that arc $l$ is incident to a node in $\widetilde{\mcl{I}}_1$.
	Since the given \ECR assignment leads to a facet-defining inequality after applying the relaxation step, its aggregation weights correspond to an extreme point of $\mcl{C}^l$ as descried in the proof of Proposition~\ref{prop:primal-weight}.
	The resulting system of equations for the associated basic feasible solution can be written as
	\begin{equation}
		\left[
		\def\arraystretch{1.2}
		\begin{array}{ccc|ccc|c}
			\pm \bar{E}_1 & \pm I_1 & \pm \bar{I}_1 \, & \, 0 & 0 & 0 \, & \, C_1 \\	
			\hline		
			0 & 0 & 0  \, & \pm \bar{E}_2 & \pm I_2 & \pm \bar{I}_2 \, & \, C_2 \\
			\hline
			\hline
			0 & 0 & 0  \, & 0 & 0 & 0 \, & \, C_3 \\
		\end{array}
		\right]
		\left[
		\begin{array}{c}
			\vc{1} \\
			\hline
			\vc{1} \\
			\hline
			\vc{0}
		\end{array}
		\right]
		=
		\left[
		\begin{array}{c}
			\pm \vc{e}^{l} \\
			\hline
			\vc{0} \\
			\hline
			\hline
			\vc{0}
		\end{array}
		\right], \label{eq:matrix-tree}
	\end{equation}
	where the columns and rows of the basis matrix have been suitably reordered. 
	In \eqref{eq:matrix-tree}, the first (resp. second) row block corresponds to bilinear terms $y_1x_i$ for arcs $i \in \mr{A}$ that are incident to the nodes in $\widetilde{\mcl{I}}_1$ (resp. $\widetilde{\mcl{I}}_2$), and the last row block corresponds to all the other bilinear terms that do not appear during aggregation. 
	The first (resp. fourth) column block denotes the transpose of the node-arc incidence matrix for nodes in $\widetilde{\mcl{I}}_1$ (resp. $\widetilde{\mcl{I}}_2$).
	The second (resp. fifth) column block contains positive or negative multiples of columns of the identity matrix representing the weights used in the relaxation step of the \ECR procedure corresponding to the McCormick bounds.
	The third (resp. sixth) column block represents positive or negative multiples of the bilinear constraints in the description of $\mcl{S}^1$ used in the relaxation step corresponding to the arcs appearing in the first (resp. second) row blocks.
	All these columns have weights equal to 1 according to Proposition~\ref{prop:primal-weight} as denoted in the first two row blocks of the solution vector multiplied with this matrix. 
	The last column in the basis corresponds to the constraints that have weights $0$ in the basic feasible solution and are added to complete the basis. 
	Lastly, $\vc{e}^{l}$ is a unit vector whose elements are all zeros except that corresponding to $y_1x_l$, which is equal to $1$.
	It follows that the linear combination of the columns in the column blocks 4, 5 and 6 of the basis matrix with weights $1$ yields the zero vector.
	This shows that these columns are linearly dependent, a contradiction.	
	\Halmos	
\endproof	

\medskip
Proposition~\ref{prop:primal-tree} implies that each non-trivial \ECR inequality can be obtained as an aggregation of constraints corresponding to a special tree structures.
The next theorem provides the converse result that aggregating constraints associated with each special tree structure can produce \ECR inequalities.

\begin{theorem} \label{thm:primal-tree-converse}
	Consider set $\mcl{S}^1$ with $\Xi$ that represents the network polytope corresponding to the network $\mr{G} = (\mr{V},\mr{A})$.
	Let $\bar{\mr{T}}$ be a tree in $\bar{\mr{G}}$ with the node set $\widetilde{\mcl{I}} \subseteq \mr{V}$, and let $l \in \mr{A}$ be an arc that is incident to exactly one node of $\bar{\mr{T}}$.
	Then, for any partition $\widetilde{\mcl{I}}_1$ and $\widetilde{\mcl{I}}_2$ of $\widetilde{\mcl{I}}$ (i.e., $\widetilde{\mcl{I}}_1 \cap \widetilde{\mcl{I}}_2 = \emptyset$ and $\widetilde{\mcl{I}}_1 \cup \widetilde{\mcl{I}}_2 = \widetilde{\mcl{I}}$), we have that
	\begin{itemize}
		\item[(i)] if $h(l) \in \widetilde{\mcl{I}}$, then $\big[\{i^+\}_{i \in \widetilde{\mcl{I}}_1},\{i^-\}_{i \in \widetilde{\mcl{I}}_2}\big]$ is an \ECR assignment for class-$l^{+}$, 
		\item[(ii)] if $h(l) \in \widetilde{\mcl{I}}$, then $\big[\{i^-\}_{i \in \widetilde{\mcl{I}}_1},\{i^+\}_{i \in \widetilde{\mcl{I}}_2}\big]$ is an \ECR assignment for class-$l^{-}$,
		\item[(iii)] if $t(l) \in \widetilde{\mcl{I}}$, then $\big[\{i^-\}_{i \in \widetilde{\mcl{I}}_1},\{i^+\}_{i \in \widetilde{\mcl{I}}_2}\big]$ is an \ECR assignment for class-$l^{+}$, 
		\item[(iv)] if $t(l) \in \widetilde{\mcl{I}}$, then $\big[\{i^+\}_{i \in \widetilde{\mcl{I}}_1},\{i^-\}_{i \in \widetilde{\mcl{I}}_2}\big]$ is an \ECR assignment for class-$l^{-}$.
	\end{itemize}
\end{theorem}	

\proof{Proof.}
	We show the result for part (i), as the proof for parts (ii)--(iv) follows from similar arguments.
	It suffices to show that the aggregation procedure performed on the constraints in the proposed assignment satisfies the \ECR conditions (C1) and (C2).
	For condition (C1), we need to show that at least $|\widetilde{\mcl{I}}_1| + |\widetilde{\mcl{I}}_2|$ bilinear terms are canceled during aggregation.
	Let $R$ be the set of arcs in $\mr{T}$, which is the directed variant of $\bar{\mr{T}}$ obtained by replacing each edge with its corresponding arc in $\mr{G}$.
	It is clear from the definition that $l \notin R$.
	As a result, for each $r \in R$, the only constraints in the aggregation that contain $x_r$ are the flow-balance constraints corresponding to the head node $h(r)$ and tail node $t(r)$ of $r$  since both of these nodes are included in $\widetilde{\mcl{I}}$ as $r$ is an arc in $\mr{T}$.
	There are four cases for the partitions of $\widetilde{\mcl{I}}$ that these head and tail nodes can belong to.
	For the first case, assume that $h(r) \in \widetilde{\mcl{I}}_1$ and $t(r) \in \widetilde{\mcl{I}}_1$.
	It follows from the \ECR assignment in case (i) that the positive flow-balance constraints $h(r)^+$ and $t(r)^+$ are used in the aggregation with weights $y_1$.
	In particular, we have $y_1\big(\sum_{k \in \delta^+(h(r))} x_k - \sum_{k \in \delta^-(h(r))\setminus \{r\}} x_k - x_r \geq f_{(h(r))}\big)$ added with $y_1\big(\sum_{k \in \delta^+(t(r))\setminus \{r\}} x_k - \sum_{k \in \delta^-(t(r))} x_k + x_r \geq f_{(t(r))}\big)$, which results in the cancellation of $y_1x_r$.
	For the second case, assume that $h(r) \in \widetilde{\mcl{I}}_1$ and $t(r) \in \widetilde{\mcl{I}}_2$.
	It follows from the \ECR assignment that the positive flow-balance constraints $h(r)^+$ and the negative flow-balance constraint $t(r)^+$ are used in the aggregation with weights $y_1$ and $(1-y_1)$, respectively.
	In particular, we have $y_1\big(\sum_{k \in \delta^+(h(r))} x_k - \sum_{k \in \delta^-(h(r))\setminus \{r\}} x_k - x_r \geq f_{(h(r))}\big)$ added with $(1-y_1)\big(-\sum_{k \in \delta^+(t(r))\setminus \{r\}} x_k + \sum_{k \in \delta^-(t(r))} x_k - x_r \geq -f_{(t(r))}\big)$, which results in the cancellation of $y_1x_r$.
	For the remaining two cases, we can use similar arguments by changing the inequality signs to show that the term $y_1x_r$ will be canceled during aggregation.
	As a result, we obtain at least $|R|$ cancellations during the aggregation corresponding to the arcs of $\mr{T}$.
	Since $\bar{\mr{T}}$ is a tree, we have that $|R| = |\widetilde{\mcl{I}}_1| + |\widetilde{\mcl{I}}_2| - 1$.
	Finally, for arc $l$, it follows from the assumption of case (i) in the problem statement that $h(l) \in \widetilde{\mcl{I}}$ and $t(l) \notin \widetilde{\mcl{I}}$.
	If $h(l) \in \widetilde{\mcl{I}}_1$, then according to the \ECR procedure for class-$l^+$, we aggregate $y_1x_l - z_l = 0$ with $y_1\big(\sum_{k \in \delta^+(h(l))} x_k - \sum_{k \in \delta^-(h(l))\setminus \{r\}} x_k - x_l \geq f_{(h(l))}\big)$, which results in the cancellation of $y_1x_l$.
	If $h(l) \in \widetilde{\mcl{I}}_2$, we aggregate $y_1x_l - z_l = 0$ with $(1-y_1)\big(-\sum_{k \in \delta^+(h(l))} x_k + \sum_{k \in \delta^-(h(l))\setminus \{r\}} x_k + x_l \geq -f_{(h(l))}\big)$, which results in the cancellation of $y_1x_l$.
	As a result, in total we have at least $|R| + 1 = |\widetilde{\mcl{I}}_1| + |\widetilde{\mcl{I}}_2|$ cancellations during the aggregation of the constraints in the \ECR assignment, showing the satisfaction of condition (C1) of the \ECR procedure. 
	For condition (C2) of the \ECR procedure, we need to show that for each constraint used in the aggregation, including the base equality, at least one bilinear term among those created after multiplication of that constraint with its corresponding weight is canceled.
	There are two types of constraints to consider.
	The first type is the flow-balance constraints in $\widetilde{\mcl{I}}_1$ and $\widetilde{\mcl{I}}_2$, which correspond to the nodes of $\bar{\mr{T}}$.
	It follows from the previous discussion that for each node $i \in \widetilde{\mcl{I}}_1 \cup \widetilde{\mcl{I}}_2$, the bilinear term $y_1x_r$, where $r$ is an arc in $\mr{T}$ that is incident to $i$, i.e., $h(r) = i$ or $t(r) = i$, is canceled during aggregation.
	This proves that at least one bilinear term is canceled in the inequality obtained after multiplying the corresponding flow-balance constraint at node $i$ with $y_1$ or $1-y_1$.
	The second type of constraints used in the aggregation is the base equality $l$.
	The proof follows from an argument similar to that given above where we showed that the bilinear term $y_1x_l$ that appears in the base constraint $y_1x_l - z_l = 0$ is canceled.
	We conclude that condition (C2) of the \ECR procedure is satisfied for all constraints used in the aggregation. 
	\Halmos
\endproof	

\medskip
In view of Theorem~\ref{thm:primal-tree-converse}, note that for the simplest choice of the tree $\bar{\mr{T}}$, i.e., an empty set, the resulting \ECR inequalities recover the classical McCormick bounds.
Therefore, considering any nonempty tree structure can potentially improve the McCormick results by adding new valid inequalities for the bilinear set.

\medskip
Proposition~\ref{prop:primal-tree} and Theorem~\ref{thm:primal-tree-converse} suggest that the \ECR assignments have a simple graphical interpretation for $\mcl{S}^1$, which can be used to generate all non-trivial \ECR inequalities to describe $\conv(\mcl{S}^1)$ without the need to search for all possible constraints and their aggregation weights that satisfy the \ECR conditions as is common for general bilinear sets.
This attribute can significantly mitigate cut-generation efforts when used systematically to produce cutting plane.
Such a systematic procedure can be designed by identifying tress of a given network and then following the result of Theorem~\ref{thm:primal-tree-converse} to obtain the corresponding \ECR assignments.
We illustrate this procedure in the following example.

\begin{figure}[!t]
	\centering
	\includegraphics[scale=3]{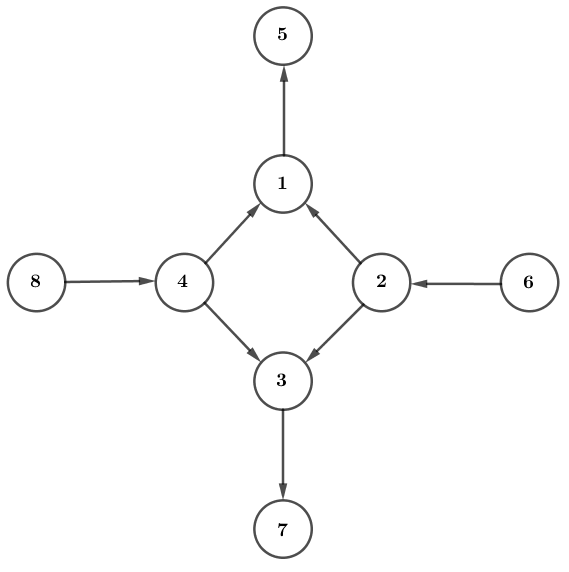}
	\caption{Graph $\mr{G}$ of Example~\ref{ex:primal-tree}}\label{fig:primal-tree}
\end{figure}

\medskip
\begin{example} \label{ex:primal-tree}
	Consider set $\mcl{S}^1$ where $\Xi$ represents the network model corresponding to a \textit{spiked cycle} graph $\mr{G} = (\mr{V},\mr{A})$ shown in Figure~\ref{fig:primal-tree}. 
	We refer to each arc in this network as a pair $(i,j)$ of its tail node $i$ and its head node $j$, and denote its corresponding flow variable as $x_{i,j}$.
	Assume that we are interested in finding \ECR assignments for class-$(1,5)^{+}$.
	According to Theorem~\ref{thm:primal-tree-converse}, we need to identify the trees that contain exactly one of the tail and head nodes of arc $(1,5)$.
	For instance, we may select the tree $\bar{\mr{T}}$ composed of the nodes $\widetilde{\mcl{I}} = \{8,4,1,2,6\}$ that contain the tail node of arc $(1,5)$.
	Consider the partitions $\widetilde{\mcl{I}}_1 = \{8,2\}$ and $\widetilde{\mcl{I}}_2 = \{4,1,6\}$.
	Following case (iii) in Theorem~\ref{thm:primal-tree-converse}, we can obtain the \ECR assignment $\big[\{8^-, 2^-\},\{4^+, 1^+, 6^+\}\big]$ for class-$(1,5)^{+}$.
	As a result, we multiply the negative flow-balance constraints at nodes $8$ and $2$ with $y_1$, and we multiply the positive flow-balance constraints at nodes $4$, $1$, and $6$ with $1-y_1$, and aggregate them with the base bilinear equality corresponding to arc $(1,5)$ with weight 1 to obtain the aggregated inequality
	\begin{multline*}
		-z_{1,5} - y_1 x_{2,3} - y_1 x_{4,3} + (f_8 + f_2 +f_1 + f_4 + f_6)y_1\\
		+ x_{1,5} -x_{2,1} + x_{4,3} - x_{8,4} + x_{6,2} - f_1 - f_4 - f_6 \geq 0
	\end{multline*}
	where $f_i$ denotes the supply/demand value at node $i$.
	Following Remark~\ref{rem:1-simplex}, we may relax each remaining bilinear term $-y_1 x_{2,3}$ and $-y_1x_{4,3}$ into three possible linear expressions, leading to 9  total \ECR inequalities.
	If implemented inside of a separation oracle, we can use Remark~\ref{rem:1-separation} to find the most violated inequality among these 9 efficiently.
\hfill	$\blacksquare$
\end{example}

\medskip

\subsubsection*{Generalization of the Coefficient Matrix for Bilinear Constraints.}
Consider a generalization of $\mcl{S}^1$ where the bilinear constraints may contain multiple bilinear terms:
\[
\widetilde{\mcl{S}}^1 = \left\{ (\vc{x};y;\vc{z}) \in \Xi \times \Delta_1 \times \Re^{\kappa} \, \middle|\,
y_1 \tilde{A}^k \vc{x} = z_{k}, \, \, \, \forall k \in K
\right\},
\]
where $\tilde{A}^k$ is a matrix of appropriate size with potentially multiple nonzero elements.
For instance, $\widetilde{\mcl{S}}^1$ may include the constraint $2y_1x_i -5 y_1x_j = z_k$ for some $i, j \in N$.
In this case, the coefficient matrix \eqref{eq:ECR-matrix} of $\mcl{C}^l$ will be modified as follows after rearranging columns and rows.
\begin{equation}
	\Bigl[
	\begin{array}{c|c|c|c|c|c|c|c}
		E^{\tr} \, & \, -E^{\tr} \, & \, I \, & \, -I \, & \, I \, & \, -I & \, \tilde{\tilde{A}} \, & \, -\tilde{\tilde{A}}
	\end{array}
	\Bigr]. \label{eq:ECR-matrix-3}
\end{equation}
In the above matrix, the row and column blocks are defined similarly to those of \eqref{eq:ECR-matrix} with a difference that the seventh and eighth column blocks correspond to the weights of the bilinear constraints $y_1\tilde{A}^k\vc{x} = z_k$, which are represented by $\beta^+_k$ and $\beta^-_k$ in the dual weight vector, respectively, for all $k \in K_l$.
In particular, the element at column $k \in K_l$ and row $i \in N$ of $\tilde{\tilde{A}}$ is equal to $\tilde{A}^k_{1i}$.
It is clear from the structure of \eqref{eq:ECR-matrix-3} that this matrix may lose the TU property when $\tilde{A}^k$ contains multiple nonzero elements for some $k \in K_l$.
In fact, this property may not hold even if $\tilde{A}^k_{1i} \in \{0,1,-1\}$ for all $k \in K_l$ and $i \in N$.
As a result, there is no guarantee that the aggregation weights for all the \ECR inequalities corresponding to non-trivial facets of $\conv(\widetilde{\mcl{S}}^1)$ will be $1$.
%Example?? shows a problem where a facet-defining \ECR inequality is obtained using aggregation weights that are not $1$. 
While an explicit derivation of the convex hull description through identifying special network structures, such as those presented for $\mcl{S}^1$, may not be attainable for this problem in its original space of variables, we can use the following ancillary result to apply Theorem~\ref{thm:primal-tree-converse} and obtain a convex hull description for $\widetilde{\mcl{S}}^1$ in a higher-dimensional space.

\begin{proposition} \label{prop:convex-high}
	Consider sets
	\[
	\mcl{S}^1 = \left\{ (\vc{x};y;\vc{w}) \in \Xi \times \Delta_1 \times \Re^{n} \, \middle|\,
	y_1x_i = w_{i}, \, \, \, \forall i \in N
	\right\},
	\] 
	and
	\[
	\mcl{D} = \left\{ (\vc{x};y;\vc{w};\vc{z}) \in \Xi \times \Delta_1 \times \Re^{n} \times \Re^{\kappa} \, \middle|\,
	\tilde{A}^k\vc{w} = z_{k}, \, \, \, \forall k \in K
	\right\},
	\] 
	Then,
	\begin{equation}
		\conv\left( (\mcl{S}^1 \times \Re^{\kappa}) \cap \mcl{D} \right) = \left(\conv(\mcl{S}^1) \times \Re^{\kappa}\right) \cap \mcl{D}. \label{eq:conv_high}
	\end{equation}
\end{proposition}	

\proof{Proof.}
	We prove the result by showing both directions of inclusion for the equality.
	The direct inclusion follows from the fact that the convex hull of intersection of two sets is a subset of the intersection of the convex hulls of those sets.
	For the reverse inclusion, we need to show that $\conv\left( (\mcl{S}^1 \times \Re^{\kappa}) \cap \mcl{D} \right) \supseteq \left(\conv(\mcl{S}^1) \times \Re^{\kappa}\right) \cap \mcl{D}$.
	Consider a point $\bar{\phi} = (\bar{\vc{x}};\bar{y};\bar{\vc{w}};\bar{\vc{z}}) \in \left(\conv\left( \mcl{S}^1 \right) \times \Re^{\kappa}\right) \cap \mcl{D}$.
	We show that $\bar{\phi} \in \conv\left( (\mcl{S}^1 \times \Re^{\kappa}) \cap \mcl{D} \right)$.
	It follows from the assumption that $\bar{z}_k = \tilde{A}^k\bar{\vc{w}}$ for all $k \in K$.
	Further, there must exist a finite collection of points $\hat{\phi}^j = (\hat{\vc{x}}^j;\hat{y}^j;\hat{\vc{w}}^j) \in \mcl{S}^1$ for $j \in J$ such that $\bar{\vc{x}} = \sum_{j \in J}\omega_j\hat{\vc{x}}^j$, $\bar{y} = \sum_{j \in J}\omega_j\hat{y}^j$, and $\bar{\vc{w}} = \sum_{j \in J}\omega_j\hat{\vc{w}}^j$ for some non-negative weights $\omega_j$ such that $\sum_{j \in J}\omega_j = 1$. 
	Consider the set of points $\dot{\phi}^j = (\dot{\vc{x}}^j;\dot{y}^j;\dot{\vc{w}}^j;\dot{\vc{z}}^j)$ for $j \in J$ such that $\dot{\vc{x}}^j = \hat{\vc{x}}^j$, $\dot{y}^j = \hat{y}^j$, $\dot{\vc{w}}^j = \hat{\vc{w}}^j$, and $\dot{\vc{z}}^j_k = \tilde{A}^k\dot{\vc{w}}^j$ for all $k \in K$.
	It is clear that $\dot{\phi}^j \in (\mcl{S}^1 \times \Re^{\kappa}) \cap \mcl{D}$ for all $j \in J$ by definition of the components of these points.
	It follows that $\bar{\phi} = \sum_{j \in J}\omega_j\dot{\phi}^j$, since $\bar{\vc{x}} = \sum_{j \in J}\omega_j\dot{\vc{x}}^j$, $\bar{y} = \sum_{j \in J}\omega_j\dot{y}^j$, and $\bar{\vc{w}} = \sum_{j \in J}\omega_j\dot{\vc{w}}^j$ by definition, and since $\bar{z}_k = \tilde{A}^k\bar{\vc{w}} = \tilde{A}^k(\sum_{j \in J}\omega_j\hat{\vc{w}}^j) = \sum_{j \in J}\omega_j \tilde{A}^k\hat{\vc{w}}^j = \sum_{j \in J}\omega_j \tilde{A}^k\dot{\vc{w}}^j = \sum_{j \in J}\omega_j \dot{\vc{z}}^j_k$ for all $k \in K$.
	This proves that $\bar{\phi} \in \conv\left( (\mcl{S}^1 \times \Re^{\kappa}) \cap \mcl{D} \right)$.
	\Halmos 
\endproof	

\medskip
The result of Proposition~\ref{prop:convex-high} shows that we can obtain a convex hull description for $\widetilde{\mcl{S}}^1$ in a higher dimension, which is expressed on the left-hand-side of \eqref{eq:conv_high}, by finding the convex hull of $\mcl{S}^1$ through application of Theorem~\ref{thm:primal-tree-converse} and then intersecting it with the linear constraints in $\mcl{D}$ as indicated on the right-hand-side of \eqref{eq:conv_high}.

%%%%%%%%%%%%%%%%%%%%%%%%%%%%%%%%%%%%%
%%%%%%%%%%%%%%%%%%%%%%%%%%%%%%%%%%%%%%%%%
%%%%%%%%%%%%%%%%%%%%%%%%%%%%%%%%%%%%%%%%%%%%%%
\subsection{The Case with $m > 1$.} \label{subsec:primal-multi}
%%%%%%%%%%%%%%%%%%%%%%%%%%%%%%%%%%%%%%%%%%%%%
%%%%%%%%%%%%%%%%%%%%%%%%%%%%%%%%%%%%%%%%%
%%%%%%%%%%%%%%%%%%%%%%%%%%%%%%%%

In this section, we consider the general case where $m > 1$ in $\mcl{S}$.
The coefficient matrix of $\mcl{C}^l$ in \eqref{eq:Cl} can be written as follows after a suitable rearrangement of columns and rows.

\begin{gather} \label{eq:ECR-matrix-2}
	\left[
	\def\arraystretch{1.5}
	\begin{array}{c|c|c|c||c||c|c||cc|cc|c|cc||cc|cc|c|cc}
		E^{\tr} \, & \vc{0} \, & \, \dotsc \, & \, \vc{0} \, & \, \, -E^{\tr} \, & I \, & \, -I \, & \, I \, & \, -I \, & \vc{0} \, & \, \vc{0} \, & \, \dotsc \, & \, \vc{0} \, & \, \vc{0} \, & \, \bar{I}^1 \, & \, -\bar{I}^1 \, & \, \vc{0} \, & \, \vc{0} \, & \, \dotsc \, & \, \vc{0} \, & \, \vc{0} \, \\
		\hline
		\, \vc{0} \, & E^{\tr} \, & \, \dotsc \, & \, \vc{0} \, & \, -E^{\tr} \, & \, I \, & \, -I \, & \vc{0} \, & \, \vc{0} \, & \, I \, & \, -I \, & \, \dotsc \, & \, \vc{0} \, & \, \vc{0} \, & \, \vc{0} \, & \, \vc{0} \, & \bar{I}^2 \, & \, -\bar{I}^2 \, & \, \dotsc \, & \, \vc{0} \, & \, \vc{0} \\
		\hline
		\, \vdots \, & \, \vdots \, & \, \ddots \, & \, \vdots \, & \, \vdots \, & \, \vdots \, & \vdots \, & \, \vdots \, & \, \vdots \, & \, \vdots \, & \, \ddots \, & \, \vdots \, & \, \vdots \, & \, \vdots \, & \, \vdots \, & \vdots \, & \, \vdots \, & \, \ddots \, & \, \vdots \, & \, \vdots \, \\
		\hline
		\, \vc{0} \, & \, \vc{0} \, & \, \dotsc \, & E^{\tr} \, & \, -E^{\tr} \, & \, I \, & \, -I \, & \vc{0} \, & \, \vc{0} \, &  \, \vc{0} \, & \, \vc{0} \, & \, \dotsc \, & \, I \, & \, -I \, & \, \vc{0} \, & \, \vc{0} \, & \, \vc{0} \, & \, \vc{0} \, & \, \dotsc \, & \bar{I}^m \, & \, -\bar{I}^m \, \\
	\end{array}
	\right]. 
\end{gather}
%In the above matrix, each row block $j \in M$ represents the bilinear terms $y_jx_i$ (i.e., $w^j_i$ in the disjunctive programming formulation \eqref{eq_ef}) for all $i \in N$.
In the above matrix, each row block $j \in M$ represents the bilinear terms $y_jx_i$ for all $i \in N$. 
The first $m$ column blocks correspond to the weights of the flow-balance constraints in $\Xi$ multiplied by $y_{j}$ for $j \in M$, which are denoted by $\gamma^j_t$ for all $t \in T$ in the dual vector.
The next column block represents the weights of the flow-balance constraints in $\Xi$ multiplied by $1-\sum_{j \in M}y_{j}$, which are denoted by $\theta_t$ for all $t \in T$ in the dual vector.
The next two column blocks indicate the lower and upper bound constraints on variables in $\Xi$ multiplied by $1-\sum_{j \in M}y_{j}$, which are recorded by $\lambda_i$ and $\mu_i$, respectively, for all $i \in N$.
The next $2m$ columns blocks correspond to the weights of the lower and upper bound constraints on variables in $\Xi$ multiplied by $y_{j}$ for $j \in M$, which are recorded by $\eta^j_i$ and $\rho^j_i$, respectively, for all $i \in N$.
The last $2m$ column blocks correspond to the weights of the positive and negative bilinear constraints in $\mcl{S}$, which are represented by $\beta^+_k$ and $\beta^-_k$ for all $k \in K_l$.
For instance, for constraint $y_jx_i - z_k \geq 0$ with $(i,j) \in N\times M$ and $k \in K_l$, the elements of column $k$ in $\bar{I}^j$ are all zero except the one in the row corresponding to the bilinear term $y_jx_i$ which is equal to one.

\medskip
It is clear from the structure of \eqref{eq:ECR-matrix-2} that this matrix does not have the TU property, implying that a result similar to that of Proposition~\ref{prop:primal-weight} does not necessarily hold.
Therefore, there is no guarantee that the aggregation weights for all the \ECR assignments obtained from the extreme points of $\mcl{C}^l$ are $1$. 
In fact, Example~\ref{ex:ECR-weight} shows that there exists \ECR assignments with aggregation weights that map to extreme points of $\mcl{C}^l$ with components that are not all $0$ or $1$. 
%As a result, a direct interpretation and explicit derivation for all \ECR inequalities based on the underlying network structures, such as trees for the case with $m=1$, may not be achievable.
%Nevertheless, we can still apply Algorithm? to obtain valid inequalities for $\conv(\mcl{S})$ considering relaxations with one of the $y$ variables at a time.

\medskip

\begin{example} \label{ex:ECR-weight}
	Consider set $\mcl{S}$ where $\Xi$ describes the network polytope corresponding to network $\mr{G} = (\mr{V},\mr{A})$ in Figure~\ref{fig:primal-tree}, and $\Delta = \{(y_1,y_2) \in \Re^2_+ | y_1 + y_2 \leq 1\}$.
	Select class-$l^+$ corresponding to the base equality $y_1x_{8,4} - z_l = 0$.
	Let $\big[\mcl{I}_1,\mcl{I}_2,\bar{\mcl{I}} \big| \mcl{J},\bar{\mcl{J}}\big]$ be an assignment for class-$l^+$ where $\mcl{I}_1 = \{4^+, 3^+\}$, $\mcl{I}_2 = \{2^-\}$, $\bar{\mcl{I}} = \{1^+\}$, $\mcl{J} = \{(4,1)\}$, and $\bar{\mcl{J}} = \{(2,3)\}$.
	Next, we show that the above assignment is an \ECR assignment for class-$l^+$ when considering the dual weight $1$ for all constraints except the bound constraint in $\mcl{J}$ which has a dual weight of $2$ in the aggregation.
	Specifically, we aggregate the base constraint $y_1x_{8,4} - z_l \geq 0$ with weight $1$, the positive flow-balance constraints at node $4$ that is $x_{4,1} + x_{4,3} - x_{8,4} \geq f_4$ with weight $y_1$, the positive flow-balance constraints at node $3$ that is $x_{3,7} - x_{4,3} - x_{2,3} \geq f_3$ with weight $y_1$, the negative flow-balance constraints at node $2$ that is $x_{6,2} - x_{2,1} - x_{2,3} \geq -f_2$ with weight $y_2$, the positive flow-balance constraints at node $1$ that is $x_{1,5} - x_{4,1} - x_{2,1} \geq f_1$ with weight $1-y_1-y_2$, the lower bound constraint for arc $(4,1)$ that is $x_{4,1} \geq 0$ with weight $2(1-y_1-y_2)$, and the upper bound constraint for arc $(2,3)$ that is $u_{2,3} - x_{2,3} \geq 0$ with weight $1-y_1-y_2$.
	During this aggregation, six bilinear terms will be canceled, satisfying condition (C1) of the \ECR procedure.
	Further, at least one bilinear term for each constraint involved in the aggregation is canceled, satisfying condition (C2) of the \ECR procedure. 
	Therefore, the above assignment is a valid \ECR assignment for class-$l^+$.
	Next, we argue that the dual weight vector for this assignment corresponds to an extreme point of $\mcl{C}^l$ in \eqref{eq:Cl}.
	For $\mcl{S}$ in this example, the coefficient matrix of $\mcl{C}^l$, as depicted in \eqref{eq:ECR-matrix-2}, has 16 rows corresponding to bilinear terms $y_jx_i$ for all $j = 1,2$ and $i \in \mr{A}$.
	It is easy to verify that the columns corresponding to the six constraints in the above \ECR assignment are linearly independent.
	As a result, we can form a basis by adding to the above six columns 10 more linearly independent columns corresponding to the constraints used in the relaxation step for the bilinear terms remaining in the aggregated inequality together with the columns that complete the basis.
	The resulting basis corresponds to a basic feasible solution where all variables (interpreted as the dual weights for constraints involved in the \ECR procedure) are $0$ or $1$, except the one associated with the column representing the lower bound constraint for arc $(4,1)$ which has a dual weight equal to $2$. 
	Therefore, there exists extreme points of $\mcl{C}^l$ with components that are not all $0$ or $1$.
	\hfill	$\blacksquare$
\end{example}

\medskip
In light of the above observation, even though identifying the aggregation weights for the \ECR procedure for $\mcl{S}$ is not as straightforward compared to that of $\mcl{S}^1$, we next show that a generalization of the previously discussed tree structures can still be detected for a given \ECR assignment.
First, we give a few definitions that will be used to derive these results.

\medskip

\begin{definition} \label{def:primal-parallel-1}
	Consider set $\mcl{S}$ where $\Xi$ describes the network polytope corresponding to network $\mr{G} = (\mr{V},\mr{A})$.
	We define a \textit{parallel network} $\mr{G}^j = (\mr{V}^j,\mr{A}^j)$ for $j \in M$ to be a replica of $\mr{G}$ that represents the multiplication of flow variables $\vc{x}$ with $y_j$ during the aggregation procedure.
	For instance, the selected nodes in the parallel network $j$ represent the multiplication of their corresponding flow-balance constraints with variable $y_j$.	
\end{definition}		

\medskip
In view of Definition~\ref{def:primal-parallel-1}, to simplify presentation, we use the same node and arc notation across all parallel networks.
In particular, given the network $\mr{G} = (\mr{V},\mr{A})$, we refer to the replica of node $v \in \mr{V}$ (resp. arc $a \in \mr{A}$) in the parallel network $\mr{G}^j$ by $v$ (resp. $a$) for all $j \in M$.
Following this rule, given a subnetwork $\dot{\mr{G}}$ of the parallel network $\mr{G}^j$ for some $j \in M$, we say that a node $v \in \mr{V}$ (resp. arc $a \in \mr{A}$) is adjacent to (resp. incident to) a node of $\mr{G}^j$, if the replica of $v$ in the parallel network $\mr{G}^j$ is adjacent to (resp. incident to) a node of $\mr{G}^j$.

%\begin{definition} \label{def:primal-parallel-2}	
%	We say that a subnetwork $\dot{\mr{G}}^1 = (\dot{\mr{V}}^1,\dot{\mr{A}}^1)$ of the parallel network $\mr{G}^i$ and a subnetwork $\dot{\mr{G}}^2 = (\dot{\mr{V}}^2,\dot{\mr{A}}^2)$ of the parallel network $\mr{G}^j$ (with possibly $i=j$) are \textit{vertically connected through a node $v \in \mr{V}$} if $v$ is adjacent to at least one node of $\dot{\mr{V}}^1$ and $v$ is adjacent to at least one node of $\dot{\mr{V}}^2$.
%	In this definition, if node $v$ is in the subnetwork, it is counted as an adjacent node to the connection node $v$.
%Therefore, a subnetwork $\dot{\mr{G}}^1$ that contains node $v$ is connected with another network through $v$ if and only if $\dot{\mr{V}}^1 = \{v\}$.
%	Similarly, we say that the two subnetwork $\dot{\mr{G}}^1$ and $\dot{\mr{G}}^2$ are \textit{vertically connected through an arc $a \in \mr{A}$} if $a$ is adjacent to at least one node of $\dot{\mr{G}}^1$ and $a$ is adjacent to at least one node of $\dot{\mr{G}}^2$.
%\end{definition}	

\medskip

\begin{definition} \label{def:primal-parallel-3}	
	Consider set $\mcl{S}$ where $\Xi$ describes the network polytope corresponding to network $\mr{G} = (\mr{V},\mr{A})$.
 Consider a collection of networks $\dot{\mr{G}}^k = (\dot{\mr{V}}^k,\dot{\mr{A}}^k)$, for $k=1,\dotsc,r$, where each $\dot{\mr{G}}^k$ is a subnetwork of a parallel networks $\mr{G}^j$ for some $j \in M$.
 We say that this collection is \textit{vertically connected through the connection nodes $C_v \subseteq \mr{V}$ and connection arcs $C_a \subseteq \mr{A}$} if there exists an ordering $s_1, s_2, \dotsc, s_r$ of indices $1,\dotsc,r$ such that for each $i=1, \dotsc, r-1$, there exists either (i) an arc $a \in C_a$ such that $a$ is incident to a node of $\dot{\mr{G}}^{s_{i+1}}$ and it is incident to a node of some subnetworks among $\dot{\mr{G}}^{s_{1}}, \dotsc, \dot{\mr{G}}^{s_{i}}$, or (ii) a set of nodes $v_1, \dotsc, v_p \in C_v$, each adjacent to the previous one in the original network $\mr{G}$, such that $v_1$ is adjacent to a node of $\dot{\mr{G}}^{s_{i+1}}$ and $v_{p}$ is adjacent to a node of some subnetworks among $\dot{\mr{G}}^{s_{1}}, \dotsc, \dot{\mr{G}}^{s_{i}}$.
	In this definition, if node $v_1$ is in $\dot{\mr{G}}^{s_{i+1}}$, it counts as an adjacent node to the connection node $v_1$.
\end{definition}	

\medskip
The next example illustrates the concepts introduced in Definitions~\ref{def:primal-parallel-1} and \ref{def:primal-parallel-3}.

\medskip
\begin{example} \label{ex:definitions}
Consider set $\mcl{S}$ where $\Xi$ describes the network polytope corresponding to network $\mr{G} = (\mr{V},\mr{A})$ in Figure~\ref{fig:primal-tree}, and $\Delta = \{(y_1,y_2) \in \Re^2_+ | y_1 + y_2 \leq 1\}$. According to Definition~\ref{def:primal-parallel-1}, since $m = 2$, we create two parallel networks $\mr{G}^j = (\mr{V}^j,\mr{A}^j)$ for $j = 1,2$. Consider subnetworks $\dot{\mr{G}}^1 = (\dot{\mr{V}}^1,\dot{\mr{A}}^1)$ and $\dot{\mr{G}}^2 = (\dot{\mr{V}}^2,\dot{\mr{A}}^2)$ of the parallel network $\mr{G}^1 = (\mr{V}^1,\mr{A}^1)$, and the subnetwork $\dot{\mr{G}}^3 = (\dot{\mr{V}}^3,\dot{\mr{A}}^3)$ of the parallel network $\mr{G}^2 = (\mr{V}^2,\mr{A}^2)$, where $\dot{\mr{V}}^1 = \{1, 5\}$, $\dot{\mr{A}}^1 = \{(1, 5)\}$, $\dot{\mr{V}}^2 = \{7\}$, $\dot{\mr{A}}^2 = \emptyset$, $\dot{\mr{V}}^3 = \{2, 6\}$, $\dot{\mr{A}}^3 = \{(6, 2)\}$. Further, define $C_v = \{3\}$ and $C_a = \{(2,1)\}$. According to Definition~\ref{def:primal-parallel-3}, the subnetworks $\dot{\mr{G}}^k$, for $k = 1,2,3$, are vertically connected through the connection nodes $C_v$ and connection arcs $C_a$ due to the following observations. Consider the ordering $s_1 = 1$, $s_2 = 3$, and $s_3 = 2$ of the indices $1, 2, 3$ of the subnetworks. It follows that $\dot{\mr{G}}^{s_2}$ is vertically connected to $\dot{\mr{G}}^{s_1}$ through arc $(2,1) \in C_a$ as it satisfies the condition (i) of Definition~\ref{def:primal-parallel-3} because it is incident to node $1$ of $\dot{\mr{G}}^{s_1}$ and node $2$ of $\dot{\mr{G}}^{s_2}$. Similarly, $\dot{\mr{G}}^{s_3}$ is vertically connected to $\dot{\mr{G}}^{s_2}$ through node $3 \in C_v$ as it satisfies the condition (ii) of Definition~\ref{def:primal-parallel-3} because it is adjacent to node $2$ of $\dot{\mr{G}}^{s_2}$ and node $7$ of $\dot{\mr{G}}^{s_3}$. 

	\hfill	$\blacksquare$	
\end{example}

\medskip
Now, we are ready to show the relationship between an \ECR assignment and a special graphical structure in the underlying network.

\begin{proposition} \label{prop:primal-forest}
	Consider set $\mcl{S}$ with $\Xi$ that represents the network polytope corresponding to the network $\mr{G} = (\mr{V},\mr{A})$.
	Let $\big[\mcl{I}_1,\dotsc,\mcl{I}_m,\bar{\mcl{I}} \big| \mcl{J},\bar{\mcl{J}}\big]$ be an \ECR assignment
	% for class-$l^{\pm}$ with $l \in K$ 
	that leads to a non-trivial facet-defining inequality of $\conv(\mcl{S})$.
	Assume that $\cup_{j \in M} \mcl{I}_j \neq \emptyset$.
	For each $j \in M$, define $\widetilde{\mcl{I}}^j = \{i \in \mr{V} | i^{\pm} \in \mcl{I}_j \}$, $\widetilde{\mcl{I}} = \{i \in \mr{V} | i^{\pm} \in \bar{\mcl{I}} \}$, and $\widetilde{\mcl{J}} = \{i \in \mr{A} | i \in \mcl{J} \cup \bar{\mcl{J}} \}$.
	%Then, in each parallel network $\mr{G}^j$ for $j \in M$, there exist a forest $\bar{\mr{F}}^j$, composed of trees $\bar{\mr{T}}^j_k$ for $k$ in an index set $\Gamma_j$, whose nodes are $\widetilde{\mcl{I}}^j$.
	%Further, the collection of all these trees are vertically connected through connection nodes $\widetilde{\mcl{I}}$ and connection arcs $\widetilde{\mcl{J}}$.
	%Lastly, the arc $k \in \mr{A}$ whose flow variable $x_k$ is multiplied with $y_j$ for some $j \in M$ in the base equality $l$ is either a connection arc or adjacent to a connection node or a node in the forest $\bar{\mr{F}}^j$. 
	Then, there exist forests $\bar{\mr{F}}^j$ in the parallel network $\mr{G}^j$ for $j \in M$, each composed of trees $\bar{\mr{T}}^j_k$ for $k \in \Gamma_j$, where $\Gamma_j$ is an index set, such that 
	\begin{itemize}
		\item[(i)] the forest $\bar{\mr{F}}^j$ is composed of the nodes in $\widetilde{\mcl{I}}^j$ for all $j \in M$,
		\item[(ii)] the collection of the trees $\bar{\mr{T}}^j_k$ for all $k \in \Gamma_j$ and all $j \in M$ are vertically connected through connection nodes $\widetilde{\mcl{I}}$ and connection arcs $\widetilde{\mcl{J}}$,
		\item[(iii)] the collection of all nodes in $\bar{\mr{F}}^j$ for all $j \in M$ together with the connection nodes $\widetilde{\mcl{I}}$ form a tree in $\mr{G}$, which has at least one incident node to each connection arc in $\widetilde{\mcl{J}}$. 
		%\item[(iv)] the arc $k \in \mr{A}$ whose flow variable $x_k$ is multiplied with $y_j$ for some $j \in M$ in the base equality $l$ is either a connection arc or adjacent to a connection node or a node in the forest $\bar{\mr{F}}^j$.
	\end{itemize}
\end{proposition}

\proof{Proof.}
	We show the result by proving conditions (i)--(iii).
	First, we may assume that the given \ECR assignment corresponds to class-$l^{\pm}$ for some $l \in K$.
	Since $\big[\mcl{I}_1,\dotsc,\mcl{I}_m,\bar{\mcl{I}} \big| \mcl{J},\bar{\mcl{J}}\big]$ leads to a non-trivial facet-defining inequality of $\conv(\mcl{S})$, its corresponding dual weights in the aggregation should represent an extreme point of $\mcl{C}^{l}$ defined in \eqref{eq:Cl}.
	This extreme point is associated with a basis in the coefficient matrix \eqref{eq:ECR-matrix-2}.
	In this basis, the subset of the column block that contains $E^{\tr}$ in the row block $j$ represents the flow-balance constraints (multiplied with $y_j$) for the nodes $i \in \widetilde{\mcl{I}}^j$, which can be viewed as the selected nodes in the parallel network $\mr{G}^j$ for $j \in M$.
	Further, the rows in the row block $j$ represent the flow variables (multiplied with $y_j$) for each arc in $\mr{G}$, which can be viewed as an arc in the parallel network $\mr{G}^j$.
	We may reorder the columns and rows of this basis corresponding to each parallel network $\mr{G}^j$ to obtain a diagonal formation composed of diagonal blocks $E^{\tr}_{j,k}$ for $k$ in an index set $\Gamma_j$. 
	It follows from this structure that the nodes corresponding to the columns of $E^{\tr}_{j,k}$ in the parallel network $\mr{G}^j$ are connected via arcs of $\mr{G}^j$ represented by the matrix rows.
	Therefore, these nodes can form a tree $\bar{\mr{T}}^j_k$ for $k \in \Gamma_j$, the collection of which represents a forest $\bar{\mr{F}}^j$ for all $j \in M$, satisfying condition (i) of the proposition statement.   
	
	\smallskip
	For condition (ii), considering the aforementioned diagonal block structure and representing the subset of column blocks of \eqref{eq:ECR-matrix-2} containing $\bar{I}^j$ and $I$ in the basis by one block with $\pm$ sign (as only one of them can be selected in the basis), we can write the resulting system of equations for the associated basic feasible solution as follows
	\begin{gather} \label{eq:ECR-matrix-4}
		\left[
		\def\arraystretch{1.2}
		\begin{array}{c|c|c|c||c||c||ccccc||ccccc||c}
			\begin{array}{ccc}
				E^{\tr}_{1,1} & 0 & 0\\
				0 & \ddots & 0 \\
				0 & 0 & E^{\tr}_{1,|\Gamma_1|}
			\end{array} & \vc{0} \, & \, \dotsc \, & \, \vc{0} \, & \, \, -E^{\tr} \, & \pm I & \, \pm I \, & \vc{0} \, & \, \dotsc \, & \, \vc{0} \, & \, \vc{0} \, & \, \pm\bar{I}^1 & \, \vc{0} \, & \, \dotsc \, & \, \vc{0} \, & \, \vc{0} \, & \, C_1 \\
			\hline
			\, \vc{0} \, & \begin{array}{ccc}
				E^{\tr}_{2,1} & 0 & 0\\
				0 & \ddots & 0 \\
				0 & 0 & E^{\tr}_{2,|\Gamma_2|}
			\end{array} & \, \dotsc \, & \, \vc{0} \, & \, -E^{\tr} \, & \, \pm I & \vc{0} \, & \, \pm I \, & \, \dotsc \, & \, \vc{0} \, & \, \vc{0} \, & \, \vc{0} \, & \pm\bar{I}^2 & \, \dotsc \, & \, \vc{0} \, & \, \vc{0} \, & \, C_2 \\
			\hline
			\, \vdots \, & \, \vdots \, & \, \ddots \, & \, \vdots \, & \, \vdots \, & \vdots \, & \, \vdots \, & \, \vdots \, & \, \ddots \, & \, \vdots & \, \vdots & \, \vdots \, & \, \vdots \, & \, \ddots \, & \, \vdots \, & \, \vdots & \, \vdots \\
			\hline
			\, \vc{0} \, & \, \vc{0} \, & \, \dotsc \, &  \begin{array}{ccc}
				E^{\tr}_{m,1} & 0 & 0\\
				0 & \ddots & 0 \\
				0 & 0 & E^{\tr}_{m,|\Gamma_m|}
			\end{array} & \, -E^{\tr} \, & \, \pm I \, & \vc{0} \, & \, \vc{0} \, & \, \dotsc \, & \, \pm I \, & \vc{0} \, & \vc{0} \, & \, \vc{0} \, & \, \dotsc \, & \pm\bar{I}^m & \, \vc{0} & \, C_m \\
			\hline
			\hline
			\, \vc{0} \, & \, \vc{0} \, & \, \dotsc \, &  \, \vc{0} \, & \, -E^{\tr} \, & \, \pm I \, & \vc{0} \, & \, \vc{0} \, & \, \dotsc \, & \vc{0} & \, \pm I \, & \vc{0} \, & \, \vc{0} \, & \, \dotsc \, & \, \vc{0} & \pm\bar{I}^{1,\dotsc,m} & \, C_{m+1} \\
			\hline
			\hline
			\, \vc{0} \, & \, \vc{0} \, & \, \dotsc \, &  \vc{0} & \, \vc{0} \, & \, \vc{0} \, & \vc{0} \, & \, \vc{0} \, & \, \dotsc \, & \, \vc{0} & \, \vc{0} \, & \vc{0} \, & \, \vc{0} \, & \, \dotsc \, & \vc{0} & \vc{0} & \, C_{m+2}
		\end{array}
		\right],
	\end{gather}	
	where the last column in the basis corresponds to the constraints that have weights $0$ in the basic feasible solution and are added to complete the basis, and where the last row represents all bilinear terms that do not appear in any constraints during aggregation.
	Further, the row block next to the last row corresponds to the bilinear terms that appear during aggregation but not in any of the selected flow-balance constraints; the matrix $\pm\bar{I}^{1,\dotsc,m}$ denotes the bilinear constraints in $\mcl{S}$ that contain these bilinear terms and could be used during the relaxation step.
	In the above basis, the column block that contains $-E^{\tr}$ represents the flow-balance constraints (multiplied with $1-\sum_{j \in M}y_j$) corresponding to the nodes in $\widetilde{\mcl{I}}$.
	Similarly, the column block that contains $\pm I$ in all rows represents the bound constraints on the flow variables (multiplied with $1-\sum_{j \in M}y_j$) corresponding to the arcs in $\widetilde{\mcl{J}}$.
	We refer to the column group composed of the columns of $E^{\tr}_{j,k}$ for any $j \in M$ and $k \in \Gamma_j$ as the column group representing the nodes of the tree $\bar{\mr{T}}^j_k$.
	It is clear from the diagonal structure of the submatrix containing $E^{\tr}_{j,k}$ that the column groups representing the nodes of $\bar{\mr{T}}^j_k$ are all \textit{arc disjoint}, i.e., there are no two columns from different groups that have a nonzero element in the same row.
	We claim that there exists an ordering $(s_1,r_1), (s_2, r_2), \dotsc, (s_h, r_h)$ of the pairs $(j,k)$ for all $j \in M$ and $k \in \Gamma_j$ such that for each column group representing the nodes of $\bar{\mr{T}}^{s_i}_{r_i}$, for all $i = 2, \dotsc, h$, there exists either (i) a column from the column blocks corresponding to $\widetilde{\mcl{J}}$ that has a nonzero element in a row corresponding to a row of $E^{\tr}_{s_i,r_i}$ and a nonzero element in a row corresponding to a row of $E^{\tr}_{s_t,r_t}$ for some $t \in \{1, \dotsc, i-1\}$, or (ii) a sequence of columns in the column block corresponding to $\widetilde{\mcl{I}}$, each sharing a nonzero element in a common row with the previous one, such that the first column has a nonzero element in a row corresponding to a row of $E^{\tr}_{s_i,r_i}$ and the last column has a nonzero element in a row corresponding to a row of $E^{\tr}_{s_t,r_t}$ for some $t \in \{1, \dotsc, i-1\}$.  
	Assume by contradiction that no such ordering exists.
	Therefore, we can partition the rows in the first $m+1$ row blocks of \eqref{eq:ECR-matrix-4} in such a way that no two rows in all columns composed of the columns in $E^{\tr}_{j,k}$ for all $j \in M$ and $k \in \Gamma_j$ and those corresponding to the columns of $\widetilde{\mcl{I}}$ and $\widetilde{\mcl{J}}$ have all their nonzero elements in the rows of one of these partitions. 
	In this case, considering that the column blocks composed of $\pm I$ and those composed of $\pm\bar{I}^j$ have exactly one nonzero element in the basis, we can compactly rewrite the system of equations for the basic feasible solution as follows:
	\begin{equation}
		\left[
		\def\arraystretch{1.2}
		\begin{array}{ccc|ccc|c} 
			\pm \bar{E}_1 & \pm I_1 & \pm\bar{I}_1 \, & \, 0 & 0 & 0 \, & \, \bar{C}_1 \\
			\hline
			0 & 0 & 0  \, & \pm \bar{E}_2 & \pm I_2 & \pm\bar{I}_2 \, & \, \bar{C}_2 \\
			\hline
			\hline
			0 & 0 & 0  \, & \, 0 & 0 & 0 \, & \, C_{m+2} \\
		\end{array}
		\right]
		\left[
		\begin{array}{c}
			\vc{+} \\
			\hline
			\vc{+} \\
			\hline
			\vc{0}
		\end{array}
		\right]
		=
		\left[
		\begin{array}{c}
			\pm \vc{e}^{l} \\
			\hline
			\vc{0} \\
			\hline
			\hline
			\vc{0}		
		\end{array}
		\right], \label{eq:matrix-forest}
	\end{equation}
	where the first and second row blocks respectively correspond to the first and second partitions discussed above. 
	In \eqref{eq:matrix-forest}, $\vc{e}^{l}$ is a unit vector whose elements are all zero except that corresponding to the row representing $y_{j'}x_{i'}$ for some $i',j'$ that satisfy $A^l_{j'i'} = 1$, which is equal to $1$.
	We may assume without the loss of generality that the row containing this nonzero element in $\vc{e}^l$ belongs to the first row block.
	All these columns except the ones in the last column block have positive weights because the associated constraints are assumed to be used in the aggregation.
	These weights are denoted by $\vc{+}$ in the first two row blocks of the vector multiplied with this matrix. 
	It follows that the linear combination of the columns in the second column block of the basis matrix with positive weights yields the zero vector.
	This shows that the columns are linearly dependent, a contradiction.
	Now, consider the ordering $(s_1,r_1), (s_2, r_2), \dotsc, (s_h, r_h)$ described above.
	It follows that for each $i = 2, \dotsc, h$, there exists either (i) a column from the column block corresponding to $\widetilde{\mcl{J}}$ that has a nonzero element in a row corresponding to a row of $E^{\tr}_{s_i,r_i}$ and a nonzero element in a row corresponding to a row of $E^{\tr}_{s_t,r_t}$ for some $t \in \{1, \dotsc, i-1\}$, or (ii) a sequence of columns in the column block corresponding to $\widetilde{\mcl{I}}$, each sharing a nonzero element in a common row with the previous one, such that the first column has a nonzero element in a row corresponding to a row of $E^{\tr}_{s_i,r_i}$ and the last column has a nonzero element in a row corresponding to a row of $E^{\tr}_{s_t,r_t}$ for some $t \in \{1, \dotsc, i-1\}$.  
	First, consider the case (i) in the above either-or argument holds for a column $k \in \widetilde{\mcl{J}}$.
	This column has nonzero elements in the rows representing arc $k$ in all subnetworks $\mr{G}^j$ for $j \in M$.
	Matrix $E^{\tr}_{s_i,r_i}$ has a nonzero element in row $k$ if the tree $\bar{\mr{T}}^{s_i}_{r_i}$ has an incident node to arc $k$.
	we conclude that $\bar{\mr{T}}^{s_i}_{r_i}$ and $\bar{\mr{T}}^{s_t}_{r_t}$ have at least one node incident to $k$, satisfying criterion (i) in Definition~\ref{def:primal-parallel-3}.
	Second, consider the case (ii) in the above either-or argument holds for a sequence $k_1, \dotsc, k_p$ of the nodes corresponding to $\widetilde{\mcl{I}}$ where $p \leq |\widetilde{\mcl{I}}|$.
	Any such column, say $k_1$, has nonzero elements in the rows representing the arcs that are incident to node $k_1$ in all parallel networks $\mr{G}^j$ for $j \in M$.
	Therefore, since each column contains a nonzero element in a common row with the previous one, the nodes corresponding to these columns must be adjacent to one another in $\mr{G}$.
	Further, since the column corresponding to $k_1$ has a nonzero element in a row corresponding to a row of $E^{\tr}_{s_i,r_i}$, we conclude that $k_1$ is adjacent to a node in $\bar{\mr{T}}^{s_i}_{r_i}$, which means that either $k_1$ belongs to this tree, or it is adjacent to a node of this tree.
	A similar argument can be made about node $k_p$ and the tree $\bar{\mr{T}}^{s_t}_{r_t}$.
	This satisfies criterion (ii) of Definition~\ref{def:primal-parallel-3}.
	This proves condition (ii) of the proposition statement due to Definition~\ref{def:primal-parallel-3}.   
	
	\smallskip
	For condition (iii), we show there exists a sequence $v_1, \dotsc, v_q$ of all the nodes in $\cup_{j \in M}\widetilde{\mcl{I}}^j \cup \widetilde{\mcl{I}}$, where $q = |\cup_{j \in M}\widetilde{\mcl{I}}^j \cup \widetilde{\mcl{I}}|$ and $v_1 \in \cup_{j \in M}\widetilde{\mcl{I}}^j$, such that every node $v_i$ is adjacent to at least one node in $v_1, \dotsc, v_{i-1}$ for every $i = 2, \dotsc, q$.
	We may assume that $v_1$ is incident to arc $i'$ defined previously that is associated with the base equality $l$.
	For other cases, where $i'$ is not incident to any nodes in $\cup_{j \in M}\widetilde{\mcl{I}}^j$, the argument will be similar with an adjustment of the partitions described below.
	It follows from the previously proven conditions (i) and (ii) of the problem statement, as well as Definition~\ref{def:primal-parallel-3} that there exists a sequence $\bar{v}_1, \dotsc, \bar{v}_p$ of all the nodes in $\cup_{j \in M}\widetilde{\mcl{I}}^j \cup \widehat{\mcl{I}}$ for some $\widehat{\mcl{I}} \subseteq \widetilde{\mcl{I}}$ such that $\bar{v}_i$ is adjacent to at least one node in $\bar{v}_1, \dotsc, \bar{v}_{i-1}$ for every $i = 2, \dotsc, p$.
	We claim that there exists a sequence $\hat{v}_1, \dotsc, \hat{v}_{q-p}$ of all the nodes in $\widetilde{\mcl{I}} \setminus \widehat{\mcl{I}}$ such that $\hat{v}_i$ is adjacent to at least one node in $\bar{v}_1, \dotsc, \bar{v}_{p}, \hat{v}_1, \dotsc, \hat{v}_{i-1}$ for every $i = 1, \dotsc, q-p$.
	Assume by contradiction that no such sequence exists.
	Then, we can use an argument similar to that of case (ii) above to partition the rows of \eqref{eq:ECR-matrix-4} in such a way that all columns corresponding to the nodes in $\bar{v}_1, \dotsc, \bar{v}_{p}, \hat{v}_1, \dotsc, \hat{v}_{t}$ have all their nonzero elements in the first partition, and the columns corresponding to all the remaining nodes $\hat{v}_{t+1}, \dotsc, \hat{v}_{q-p}$ have all their nonzero elements in the second partition.
	Then, a similar argument to that following \eqref{eq:matrix-forest} will show that the columns in the second group are linearly dependent, a contradiction.
	The case for the second part of the statement regarding the connection arcs in $\widetilde{\mcl{J}}$ can be shown similarly.	   
	\Halmos
\endproof

\medskip
%Even though the \ECR procedure provides a systematic tool to obtain valid inequalities for $\conv(\mcl{S})$, it requires the selection of the constraints together with their weights to be used in the aggregation.
%This task can be computationally cumbersome for general sets as it requires searching for weights that satisfy the \ECR conditions.
Although Proposition~\ref{prop:primal-forest} shows that an \ECR assignment corresponds to a forest structure in the parallel networks created from the underlying network in $\mcl{S}$, the converse result---similar to the one presented for set $\mcl{S}^1$---does not hold here.
More specifically, a forest structure that satisfies the conditions of Proposition~\ref{prop:primal-forest} does not necessarily lead to a valid \ECR assignment, and
even when it does, the calculation of the aggregation weights to satisfy the \ECR conditions is not as straightforward; see Example~\ref{ex:forest} below.

\medskip
\begin{example} \label{ex:forest}
	Consider set $\mcl{S}$ where $\Xi$ describes the network polytope corresponding to network $\mr{G} = (\mr{V},\mr{A})$ in Figure~\ref{fig:primal-tree}, and $\Delta = \{(y_1,y_2) \in \Re^2_+ | y_1 + y_2 \leq 1\}$.
	Select class-$l^+$ corresponding to the base equality $y_1x_{1,5} - z_l = 0$.
	In the parallel network $\mr{G}^1$, we select the forest $\bar{\mr{F}}^1$ composed of the tree $\bar{\mr{T}}^1_1$ with the node set $\widetilde{\mcl{I}}^1 = \{1, 2\}$.
	In the parallel network $\mr{G}^2$, we select the forest $\bar{\mr{F}}^2$ composed of the tree $\bar{\mr{T}}^2_1$ with the node set $\widetilde{\mcl{I}}^2 = \{2, 6\}$.
	We select the connection node set $\widetilde{\mcl{I}} = \{6\}$, and the connection arc set $\widetilde{\mcl{J}} = \emptyset$.
	It is easy to verify that these sets satisfy the conditions (i)--(iii) of Propositions~\ref{prop:primal-forest}.
	However, we cannot find an aggregation weight for the flow-balance constraints corresponding to the nodes in the above sets that yields a cancellation of at least 5 bilinear terms.
	As a result, there is no \ECR assignment that matches the considered forest structure.
\hfill	$\blacksquare$	
\end{example}

\medskip
A common way to circumvent the above-mentioned difficulty in obtaining valid \ECR assignments and their aggregation weights is to aim at a special class of \ECR assignments with more specific attributes that can be used to strengthen the connection between an \ECR assignment and its corresponding network structure.
An important example of such class is the class of \ECR assignments that are obtained through \textit{pairwise cancellation}.
%In this procedure, each cancellation of bilinear terms is obtained by aggregating two constraints, and each remaining bilinear term appears in only one constraint during the aggregation.
In this procedure, each cancellation of bilinear terms is obtained by aggregating two constraints.
This definition includes the bilinear terms that are \textit{canceled} during the relaxation step, i.e., the constraint used to relax the remaining bilinear terms counts as one of the two constraints in the preceding statement. 
Following this procedure, the aggregation weight for each constraint can be determined successively as the constraint is added to the assignment to ensure the satisfaction of the \ECR conditions.
The next result shows that the aggregation weights for all constraints used in the \ECR assignments obtained through pairwise cancellation are $1$.

\begin{proposition} \label{prop:pairwise ECR}
	Consider set $\mcl{S}$ where $\Xi$ describes the network polytope corresponding to network $\mr{G} = (\mr{V},\mr{A})$.
	Let $\big[\mcl{I}_1,\dotsc,\mcl{I}_m,\bar{\mcl{I}} \big| \mcl{J},\bar{\mcl{J}}\big]$ be an \ECR assignment for class-$l^{\pm}$ for some $l \in K$ corresponding to pairwise cancellation.
	Then, the aggregation weights for all constraints used in this assignment are $1$. 
\end{proposition}	

\proof{Proof.}
	Let $\bar{\vc{\pi}}^l$ be the solution vector for the system of equations \eqref{eq:Cl} of $\mcl{C}^l$ corresponding to the aggregation weights of the given \ECR assignment. 
	We may rewrite this system of equations as follows by rearranging its rows and columns.
	\begin{equation}
		\left[
		\def\arraystretch{1.5}
		\begin{array}{c|c|c} 
			P_1 \, & \, 0 \, & \, C_1 \\
			\hline
			P_2 \, & \, \pm I \, & \, C_2 \\
			\hline
			0 \, & \, 0 \, & \, C_3 \\
		\end{array}
		\right]
		\left[
		\def\arraystretch{1.5}
		\begin{array}{c}
			\bar{\vc{\pi}}^l_1 \\
			\hline
			\bar{\vc{\pi}}^l_2 \\
			\hline
			\vc{0}
		\end{array}
		\right]
		=
		\left[
		\def\arraystretch{1.5}
		\begin{array}{c}
			\pm \vc{e}^{l} \\
			\hline
			\vc{0} \\
			\hline
			\vc{0} 
		\end{array}
		\right], \label{eq:matrix-pairwise}
	\end{equation}
	In the coefficient matrix of \eqref{eq:matrix-pairwise}, the first row block represents the bilinear terms that are canceled during aggregation.
	The second row block corresponds to the remaining bilinear terms in the aggregated inequality that are relaxed in the last step of the \ECR procedure.
	The last row block represents all the bilinear terms that are not involved in the aggregation procedure.
	Further, in this matrix, the first column block corresponds to the constraints used in the aggregation, whose aggregation weights in the solution vector $\bar{\vc{\pi}}^l$ are denoted by $\bar{\vc{\pi}}^l_1$.
	The second column block corresponds to the variable bound constraints in $\Xi$ as well as the bilinear constraints in $\mcl{S}$ used in the \ECR procedure to relax the remaining bilinear terms in the aggregated inequality, whose weights in the solution vector $\bar{\vc{\pi}}^l$ are denoted by $\bar{\vc{\pi}}^l_2$.
	The last column block represents all other constraints that are not used during the \ECR procedure and their weights in the solution vector $\bar{\vc{\pi}}^l$ are zero. 
	Finally, $\vc{e}^{l}$ on the right-hand-side of this system is a unit vector whose elements are all zeros except that corresponding to the row representing $y_{j'}x_{i'}$ for some $i',j'$ that satisfy $A^l_{j'i'} = 1$, which is equal to $1$.
	It is clear that this row belongs to the first row block since according to the \ECR condition (C2), the bilinear term in the base equality $l$ must be canceled during the aggregation procedure when the assignments are not empty.
	It follows from the equation \eqref{eq:matrix-pairwise} that $P_1 \bar{\vc{\pi}}^l_1 = \pm \vc{e}^{l}$.
	Next, we analyze the structure of $P_1$.
	Note that all elements of $P_1$ belong to $\{0, -1, 1\}$ because it is a submatrix of \eqref{eq:ECR-matrix-2} that represents the coefficients of the constraints in $\mcl{S}$.
	Considering that the columns of $P_1$ represent the constraints used in the aggregation except the base equality (as that constraint has been moved to the right-hand-side to form $\mcl{C}^l$), and that the rows of $P_1$ correspond to the canceled bilinear terms during aggregation, according to condition (C1) of \ECR, we conclude that the number of rows of $P_1$ is no smaller that the number of columns of $P_1$.
	Further, it follows from condition (C2) of \ECR that each constraint used in the aggregation (after being multiplied with its corresponding weight) will have at least one bilinear term canceled, which implies that each column of $P_1$ has at least one nonzero element.
	The assumption of pairwise cancellation for the given \ECR assignment implies that each canceled bilinear term corresponding to the rows of $P_1$ are obtained through aggregation of exactly two constraints.
	As a result, each row of $P_1$ must contain exactly two nonzero elements, except for the row corresponding to the bilinear term $y_{j'}x_{i'}$ that appears in the base equality $y_{j'}x_{i'} - z_l = 0$ which must have only one nonzero element because the weight of the base equality has been fixed at $\pm 1$ and its column has been moved to the right-hand-side of the equation captured by $\pm \vc{e}^{l}$; see the derivation of \eqref{eq:Cl}.
	Therefore, we may rearrange the rows and columns of the matrices in this equation to obtain the form:
	\begin{equation}
		\left[
		\def\arraystretch{1.2}
		\begin{array}{c|c|c|c} 
			\pm 1 & 0 & \cdots & 0 \\
			\hline
			\{0, \pm 1\} & \pm 1 & \cdots & 0 \\
			\hline
			\vdots & \vdots & \ddots & \vdots \\
			\hline
			\{0, \pm 1\} & \{0, \pm 1\} & \cdots & \pm 1 \\
			\hline
			\{0, \pm 1\} & \{0, \pm 1\} & \cdots & \{0, \pm 1\} \\
			\hline
			\vdots & \vdots & \ddots & \vdots \\
			\hline
			\{0, \pm 1\} & \{0, \pm 1\} & \cdots & \{0, \pm 1\}
		\end{array}
		\right]
		\bar{\bar{\vc{\pi}}}^l_1 =	\pm \vc{e}^{1}, \label{eq:matrix-pairwise-2}
	\end{equation}
	where $\bar{\bar{\vc{\pi}}}^l_1$ is composed of the elements of $\bar{\vc{\pi}}^l_1$ that are rearranged to match the rearrangement of columns of $P_1$ in the above form, and where the first row corresponds to the bilinear term  $y_{j'}x_{i'}$ so that the right-hand-side vector becomes $\pm \vc{e}^1$.
	It follows from the above discussion about the structure of $P_1$ and equation \eqref{eq:matrix-pairwise-2} that all components of $\bar{\bar{\vc{\pi}}}^l_1$ must be equal to 1 as they need to be nonnegative.
	Finally, for the equations in the second row block of \eqref{eq:matrix-pairwise}, we have that $P_2 \bar{\vc{\pi}}^l_1 \pm I \bar{\vc{\pi}}^l_2 = \vc{0}$.
	It follows from the pairwise cancellation assumption that each row of $P_2$ contains exactly one nonzero element as it corresponds to a remaining bilinear term in the aggregation inequality.
	Since all of the elements in $P_2$ belong to $\{0, -1, 1\}$, and all the components in $\bar{\vc{\pi}}^l_1$ are equal to 1, it must hold that $\bar{\vc{\pi}}^l_2 = \vc{1}$.
	\Halmos
\endproof	

\medskip
\begin{remark} \label{rem:pairwise-primal}
	Consider the special case with $m=1$. As described in the proof of Proposition~\ref{prop:primal-tree}, for any $i \in N \setminus \{l\}$ (resp. for $i = l$), there are two (resp. three) possible scenarios for constraints that could be used in the aggregation to cancel a bilinear term $y_1x_i$.
	Since the aggregation weights for all constraints are $1$ in this case (see Proposition~\ref{prop:primal-weight}), we conclude that each cancellation is obtained through aggregation of exactly two constraints.
	Further, any remaining bilinear term in the aggregated inequality corresponds to an arc that is incident to exactly one node of the tree associated with the \ECR assignment (see Proposition~\ref{prop:primal-tree}), which implies that each such bilinear term appears in exactly one constraint during aggregation.
	As a result, \textit{all} \ECR inequalities for the case with $m = 1$ can be obtained through pairwise cancellation.
\end{remark}		

\medskip
Although the \ECR inequalities obtained through pairwise cancellation do not necessarily produce a full convex hull description for $\mcl{S}$, the result of Proposition~\ref{prop:pairwise ECR} provides three important advantages: (i) it generalizes the convexification results for the case with $m = 1$ as described in Remark~\ref{rem:pairwise-primal}; (ii) it can produce inequalities stronger than those obtained by applying Theorem~\ref{thm:primal-tree-converse} to relaxations of $\mcl{S}$ that contain one $y$ variable at a time, because it considers all the $y$ variables in their original simplex set $\Delta_m$; and (iii) it enables us to derive explicit \ECR inequalities cognizant of the underlying network structure without the need to search for the aggregation weights that satisfy the \ECR conditions, as will be shown in the sequel.
These advantages are corroborated by the computational experiments presented in Section~\ref{sec:computation}. 
The next proposition shows that the pairwise cancellation property provides more information about the forest structure presented in Proposition~\ref{prop:primal-forest}.

\begin{proposition} \label{prop:primal-forest-pairwise}
	Consider the setting of Proposition~\ref{prop:primal-forest}, and assume that the \ECR assignment $\big[\mcl{I}_1,\dotsc,\mcl{I}_m,\bar{\mcl{I}} \big| \mcl{J},\bar{\mcl{J}}\big]$ has the pairwise cancellation property.
	Further, let this assignment correspond to a class-$l^{\pm}$ for some $l \in K$ such that $A^l_{j'i'} = 1$ for some $(i',j') \in N \times M$.
	Then, in addition to the outcome of Proposition~\ref{prop:primal-forest}, we have that
	\begin{itemize}
		\item[(i)] arc $i'$ is either in $\widetilde{\mcl{J}}$ or incident to exactly one node in $\widetilde{\mcl{I}} \cup \widetilde{\mcl{I}}^{j'}$, but not both,
		\item[(ii)] each arc in $\widetilde{\mcl{J}}$ is incident to at most one node in $\widetilde{\mcl{I}} \cup \widetilde{\mcl{I}}^{j}$ for each $j \in M$,
		\item[(iii)] each node in $\widetilde{\mcl{I}} \cap \widetilde{\mcl{I}}^j$, for $j \in M \setminus \{j'\}$ (resp. $j = j'$), is adjacent to no other nodes in that set and no arcs in $\widetilde{\mcl{J}}$ (resp. $\widetilde{\mcl{J}} \cup \{i'\}$).
	\end{itemize}
\end{proposition}

\proof{Proof.}
	For case (i), it follows from condition (C2) of the \ECR procedure that the bilinear term $y_{j'}x_{i'}$ in the base equality must be canceled during aggregation.
	Further, according to the pairwise cancellation property, there must be exactly one constraint in the aggregation in addition to the base equality that would contain $y_{j'}x_{i'}$ after multiplication with the corresponding dual weight.
	There are two possible scenarios.
	The first possibility is that the bound constraints for $x_{i'}$ are used in the aggregation, which implies that arc $i'$ is a connection arc and belongs to $\widetilde{\mcl{J}}$.
	The second possibility is that the flow-balance constraint at either node $t(i')$ or $h(i')$, but not both, is used in the aggregation, which implies that $i'$ is incident to exactly one node in $\widetilde{\mcl{I}} \cup \widetilde{\mcl{I}}^{j'}$.
	
	\smallskip
	For case (ii), consider an arc $i \in \widetilde{\mcl{J}}$.
	Therefore, either of the bound constraints $x_i \geq 0$ or $u_i-x_i \geq 0$ is used in the aggregation with weight $1-\sum_{j \in M} y_j$.
	It follows from the pairwise cancellation property that, for each $j \in M$, there can be at most one additional constraint in the aggregation that contains a term $y_j x_i$.
	The only possibility for such a constraint is the flow-balance constraint at either node $t(i')$ or $h(i')$, but not both.
	We conclude that $i$ is incident to at most one node in $\widetilde{\mcl{I}} \cup \widetilde{\mcl{I}}^{j}$ for each $j \in M$. 
	
	\smallskip
	For case (iii), consider a node $i \in \widetilde{\mcl{I}} \cap \widetilde{\mcl{I}}^j$ for some $j \in M \setminus \{j'\}$.
	Therefore, the aggregation contains the (positive or negative) flow-balance constraint at node $i$ multiplied with $1-\sum_{j \in M} y_j$ due to $i \in \widetilde{\mcl{I}}$, together with the (positive or negative) flow-balance constraint at node $i$ multiplied with $y_j$ due to $i \in \widetilde{\mcl{I}}^j$.
	Therefore, the bilinear terms $y_j x_k$ for all $k \in \delta^+(i) \cup \delta^-(i)$ already appear in two constraints, which implies that they cannot appear in any other constraints during aggregation. 
	As a result, the bound constraints for each variable $x_k$ corresponding to arc $k$ cannot be included in $\widetilde{\mcl{J}}$.
	Similarly, the flow-balance constraint at node $h(k)$ for any $k \in \delta^+(i)$ and at node $t(k)$ for any $k \in \delta^-(i)$ cannot be included in the aggregation, which implies that $i$ cannot be adjacent to any other nodes in $\widetilde{\mcl{I}} \cap \widetilde{\mcl{I}}^j$.
	The proof for the case where $j = j'$ follows from a similar argument.
	\Halmos
\endproof	 

\medskip
As noted earlier, an important consequence of the pairwise cancellation property is providing the ability to derive the converse statement to those of Proposition~\ref{prop:primal-forest} and \ref{prop:primal-forest-pairwise}, which identifies \ECR assignments based on a special forest structure in the underlying network.
The procedure to obtain such an \ECR assignment is given in Theorem~\ref{thm:primal-forest-converse}, which makes use of Algorithm~\ref{alg:primal-ECR} to determine the aggregation weights for deriving the corresponding \ECR inequalities.

\begin{theorem} \label{thm:primal-forest-converse}
	Consider set $\mcl{S}$ with $\Xi$ that represents the network polytope corresponding to the network $\mr{G} = (\mr{V},\mr{A})$.
	Let $\bar{\mr{F}}^j$, for each $j \in M$, be a forest in the parallel network $\mr{G}^j$, composed of trees $\bar{\mr{T}}^j_k$ for $k \in \Gamma_j$, where $\Gamma_j$ is an index set, that satisfies the conditions (i)--(iii) of Propositions~\ref{prop:primal-forest} and \ref{prop:primal-forest-pairwise} with the corresponding node sets $\widetilde{\mcl{I}}^j$, the connection node set $\widetilde{\mcl{I}}$, the connection arc set $\widetilde{\mcl{J}}$, and the class $l$.
	Then, the assignment $\big[\mcl{I}_1,\dotsc,\mcl{I}_m,\bar{\mcl{I}} \big| \mcl{J},\bar{\mcl{J}}\big]$ obtained from Algorithm~\ref{alg:primal-ECR} is an \ECR assignment for class-$l^{\pm}$.
\end{theorem}

\proof{Proof.}
 First, we argue that conditions (i)--(iii) of Proposition~\ref{prop:primal-forest} imply that each member of the sets $\widetilde{\mcl{I}}$, $\widetilde{\mcl{J}}$, and $\widetilde{\mcl{I}}^j$ for $j \in M$ receives a label assignment through the steps of Algorithm~\ref{alg:primal-ECR}, i.e., the member is added to set $\mt{D}$ defined in that algorithm.
 It follows from condition (i) of Proposition~\ref{prop:primal-forest} that once a member of the node subset in $\widetilde{\mcl{I}}^j$, for each $j \in M$, that represents a tree $\bar{\mr{T}}^j_k$ with $k \in \Gamma_j$ is added to $\mt{D}$, all the remaining nodes in $\bar{\mr{T}}^j_k$ are eventually added to $\mt{D}$ because of the loop in lines 11--13 in the algorithm, as all nodes of the tree are connected.
 Condition (ii) of Proposition~\ref{prop:primal-forest} implies that all trees $\bar{\mr{T}}^j_k$ for $k \in \Gamma_j$ and $j \in M$ are connected through an appropriate sequence of the tree nodes, the connection nodes in $\widetilde{\mcl{I}}$, and the connection arcs in $\widetilde{\mcl{J}}$.
 Consequently, the loops in lines 10--44 of the algorithm ensure that each member of these sets is visited following that sequence and becomes added to $\mt{D}$.
 Further, condition (iii) of Proposition~\ref{prop:primal-forest} suggests that each member in the sets $\widetilde{\mcl{I}}$ and $\widetilde{\mcl{J}}$ is connected to the subgraph composed of the set of all tree nodes in $\widetilde{\mcl{I}}^j$ and their associated connection nodes and connection arcs.
 As a result, there exists a sequence of adjacent nodes and arcs that lead to each member of $\widetilde{\mcl{I}}$ and $\widetilde{\mcl{J}}$, thereby getting added to $\mt{D}$.
 
  \smallskip 
 Second, we show that each bilinear term created during the aggregation can appear in at most two constraints.
 There are four cases.
 In case 1, consider the bilinear term $y_{j'}x_{i'}$ that appears in the base equality $l$.
 Condition (i) of Proposition~\ref{prop:primal-forest-pairwise} implies that this bilinear term can appear in exactly one other constraint, which could be either the bound constraint on variable $x_{i'}$ (which would be included in $\widetilde{\mcl{J}}$) or the flow-balance constraint at one of the incident nodes to $i'$ (which would be included in $\widetilde{\mcl{I}}^{j'} \cup \widetilde{\mcl{I}}$).
 In case 2, consider a bilinear term $y_j x_i$, for some $j \in M$, that appears in the bound constraint on variable $x_i$ for any arc $i \in \widetilde{\mcl{J}}$.
 Condition (ii) of Proposition~\ref{prop:primal-forest-pairwise} implies that this bilinear term can appear in at most one other constraint, which could be the flow-balance constraint at one of the incident nodes to $i$ (which would be included in $\widetilde{\mcl{I}}^{j} \cup \widetilde{\mcl{I}}$).
 In case 3, consider a bilinear term $y_j x_i$, for some $j \in M$, that appears in the flow-balance constraint at an incident node of arc $i$ after being multiplied with both $y_j$ (i.e., the node being in $\widetilde{\mcl{I}}^j$) and $1-\sum_{j \in M}y_j$ (i.e., the node being in $\widetilde{\mcl{I}}$).
 Condition (iii) of Proposition~\ref{prop:primal-forest-pairwise} implies that this bilinear term cannot appear in any other constraints during aggregation.
 In case 4, consider a bilinear term $y_j x_i$, for some $j \in M$, that appears in the flow-balance constraint at an incident node of arc $i$ that is not in $\widetilde{\mcl{I}}^j \cap \widetilde{\mcl{I}}$.
 It follows from condition (iii) of Proposition~\ref{prop:primal-forest} that this bilinear term can appear in at most one other constraint because of the tree structure of all the nodes in $\cup_{j \in M} \widetilde{\mcl{I}}^j \cup \widetilde{\mcl{I}}$.
 
  \smallskip 
 Third, we discuss that, for any $k \in \mt{D}$ that has been newly added to this set, its label value has been determined through lines 10--44 of Algorithm~\ref{alg:primal-ECR} in such a way that, for a member $i \in \cup_{j \in M} \widetilde{\mcl{I}}^j \cup \widetilde{\mcl{I}} \cup \widetilde{\mcl{J}}$ that has been previously added to $\mt{D}$ and is adjacent/incident to $i$, the bilinear term that commonly appears in the weighted constraints corresponding to both $i$ and $k$ is canceled. 
 For instance, consider the case where $i \in \widetilde{\mcl{I}}$ (line 22 of the algorithm) and $k \in \widetilde{\mcl{I}}^j$ for some $j \in M$ is an adjacent node to $i$ (line 26 of the algorithm).
 Assume that $\mt{l}(i) = +$, and that arc $a \in \mr{A}$ is such that $t(a) = i$ and $h(a) = k$.
 It follows from line 27 of the algorithm that $\mt{l}(k) = -$.
 Considering the assignment rule in lines 48 and 52 of the algorithm, we should aggregate the constraint $\sum_{r \in \delta^+(i) \setminus \{p\}} x_r - \sum_{r \in \delta^-(i)} x_r + x_p \geq f_i$ with weight $1 - \sum_{j \in M} y_j$, together with the constraint $-\sum_{r \in \delta^+(k)} x_r + \sum_{r \in \delta^-(k) \setminus \{p\}} x_r + x_p \geq - f_k$ with weight $y_j$, which results in the cancellation of the bilinear term $y_j x_p$.
 A similar argument can be made for any other possible case in Algorithm~\ref{alg:primal-ECR}.
  
 \smallskip 
 Combining all the results shown in the previous parts, i.e., (I) each member of the sets $\widetilde{\mcl{I}}$, $\widetilde{\mcl{J}}$, and $\widetilde{\mcl{I}}^j$ for $j \in M$ receives a label assignment and is added to $\mt{D}$; (II) each bilinear term created during the aggregation can appear in at most two constraints; and (III) for any $k \in \mt{D}$, its label value is determined in such a way that the bilinear term that is common between the weighted constraints corresponding to $i$ and a previously added member $k$ in $\mt{D}$ is canceled, we conclude that at least $|\widetilde{\mcl{I}}| +  |\widetilde{\mcl{J}}|+\sum_{j \in M}|\widetilde{\mcl{I}}^j|$ bilinear terms will be canceled during aggregation in the desired assignment $\big[\mcl{I}_1,\dotsc,\mcl{I}_m,\bar{\mcl{I}} \big| \mcl{J},\bar{\mcl{J}}\big]$.
 This satisfies the \ECR conditions (C1).
 Finally, the above argument also implies that each flow-balance constraint at the nodes in $\cup_{j \in M} \widetilde{\mcl{I}}^j \cup \widetilde{\mcl{I}}$, and each variable bound constraint for the arcs in $\widetilde{\mcl{I}}$ will have at least one of their bilinear terms (after being multiplied with appropriate weights) canceled because each such node or arc will eventually be added to $\mt{D}$ when it receives a label for the desired cancellation. 
 This satisfies the \ECR condition (C2).
 We conclude that $\big[\mcl{I}_1,\dotsc,\mcl{I}_m,\bar{\mcl{I}} \big| \mcl{J},\bar{\mcl{J}}\big]$ is an \ECR assignment.
	\Halmos
\endproof

\begin{algorithm} % enter the algorithm environment
	\scriptsize
	\caption{Derive an \ECR assignment associated with a forest structure} % give the algorithm a caption
	\label{alg:primal-ECR} % and a label for \ref{} commands later in the document
	\begin{algorithmic}[1] % enter the algorithmic environment
		\REQUIRE network $\mr{G} = (\mr{V},\mr{A})$, forest node sets $\widetilde{\mcl{I}}^j$ for $j \in M$, connection node set $\widetilde{\mcl{I}}$, connection arc set $\widetilde{\mcl{J}}$, class $l$ with $(i',j') \in N \times M$ such that $A^{l}_{j'i'} = 1$, and a class sign indicator $\pm$, that satisfy the conditions of Theorem~\ref{thm:primal-forest-converse} 
		\ENSURE the \ECR assignment $\big[\mcl{I}_1,\dotsc,\mcl{I}_m,\bar{\mcl{I}} \big| \mcl{J},\bar{\mcl{J}}\big]$ for class-$l^{\pm}$
		\STATE assign an empty label denoted by $\mt{l}(i)$ to each node $i \in \widetilde{\mcl{I}} \cup_{j \in M}\widetilde{\mcl{I}}^j$ and each arc $i \in \widetilde{\mcl{J}} \cup \{l\}$, and define set $\mt{D} = \emptyset$
		\STATE set $\mathtt{l}(i') = $ the class sign indicator, and let $k \in \widetilde{\mcl{I}}^{j'} \cup \widetilde{\mcl{I}}$ either be an incident node to $i'$ or be the arc $i' \in \widetilde{\mcl{J}}$ 
		\STATE \textbf{if} $\left(k = h(i') \in \widetilde{\mcl{I}}^{j'}\right)$ or $\left(k = t(i') \in \widetilde{\mcl{I}}\right)$ or $\left(k = i' \in \widetilde{\mcl{J}}\right)$	\textbf{then}	
			\STATE \quad set $\mt{l}(k) = \mt{l}(i')$, and add $k$ to $\mt{D}$
		\STATE \textbf{else if} $\left(k = t(i') \in \widetilde{\mcl{I}}^{j'}\right)$ or $\left(k = h(i') \in \widetilde{\mcl{I}}\right)$ \textbf{then}
			\STATE \quad set $\mt{l}(k) = \neg \mt{l}(i')$, and add $k$ to $\mt{D}$ ($\neg$ represents the negation symbol)
		\STATE \textbf{end if}
		\WHILE{$\mt{D} \neq \emptyset$}
			\STATE \quad select $i \in \mt{D}$
			\STATE \quad \textbf{if} $i \in \widetilde{\mcl{I}}^j$ for some $j \in M$ \textbf{then}
				\STATE \quad \quad \textbf{for} each unlabeled node $k \in \widetilde{\mcl{I}}^j$ and $\bar{k} \in \widetilde{\mcl{I}}$ that is adjacent to $i$ \textbf{do}
					\STATE \quad \quad \quad set $\mt{l}(k) = \mt{l}(i)$, set $\mt{l}(\bar{k}) = \neg \mt{l}(i)$, and add $k$ and $\bar{k}$ to $\mt{D}$
				\STATE \quad \quad \textbf{end for}
				\STATE \quad \quad \textbf{for} each unlabeled arc $k \in \widetilde{\mcl{J}}$ that is incident to $i$ \textbf{do}
					\STATE \quad \quad \quad \textbf{if} $i = t(k)$ \textbf{then}
					\STATE \quad \quad \quad \quad set $\mt{l}(k) = \mt{l}(i)$, and add $k$ to $\mt{D}$
					\STATE \quad \quad \quad \textbf{else}
					\STATE \quad \quad \quad \quad set $\mt{l}(k) = \neg \mt{l}(i)$, and add $k$ to $\mt{D}$
					\STATE \quad \quad \quad \textbf{end if}
				\STATE \quad \quad \textbf{end for}
			\STATE \quad \textbf{end if} 
			\STATE \quad \textbf{if} $i \in \widetilde{\mcl{I}}$ \textbf{then}
				\STATE \quad \quad \textbf{for} each unlabeled node $k \in \widetilde{\mcl{I}}^j$ for $j \in M$ such that $k = i$ \textbf{do}
					\STATE \quad \quad \quad set $\mt{l}(k) = \mt{l}(i)$, and add $k$ to $\mt{D}$
				\STATE \quad \quad \textbf{end for}
				\STATE \quad \quad \textbf{for} each unlabeled node $k \in \widetilde{\mcl{I}}^j$ for $j \in M$, and $\bar{k} \in \widetilde{\mcl{I}}$ that is adjacent to $i$ \textbf{do}
					\STATE \quad \quad \quad set $\mt{l}(k) = \neg \mt{l}(i)$, set $\mt{l}(\bar{k}) = \mt{l}(i)$, and add $k$ and $\bar{k}$ to $\mt{D}$
				\STATE \quad \quad \textbf{end for}
				\STATE \quad \quad \textbf{for} each unlabeled arc $k \in \widetilde{\mcl{J}}$ that is incident to $i$ \textbf{do}
				\STATE \quad \quad \quad \textbf{if} $i = t(k)$ \textbf{then}
					\STATE \quad \quad \quad \quad set $\mt{l}(k) = \neg \mt{l}(i)$, and add $k$ to $\mt{D}$
				\STATE \quad \quad \quad \textbf{else}
					\STATE \quad \quad \quad \quad set $\mt{l}(k) = \mt{l}(i)$, and add $k$ to $\mt{D}$
				\STATE \quad \quad \quad \textbf{end if}
				\STATE \quad \quad \textbf{end for}
			\STATE \quad \textbf{end if} 
			\STATE \quad \textbf{if} $i \in \widetilde{\mcl{J}}$ \textbf{then}
				\STATE \quad \quad \textbf{for} each unlabeled node $k = t(i) \in \widetilde{\mcl{I}}^j$ and $\bar{k} = h(i) \in \widetilde{\mcl{I}}^j$ for $j \in M$ \textbf{do}
					\STATE \quad \quad \quad set $\mt{l}(k) = \mt{l}(i)$, set $\mt{l}(\bar{k}) = \neg \mt{l}(i)$, and add $k$ and $\bar{k}$ to $\mt{D}$
				\STATE \quad \quad \textbf{end for}
				\STATE \quad \quad \textbf{for} each unlabeled node $k = t(i) \in \widetilde{\mcl{I}}$ and $\bar{k} = h(i) \in \widetilde{\mcl{I}}$ \textbf{do}
					\STATE \quad \quad \quad set $\mt{l}(k) = \neg \mt{l}(i)$, set $\mt{l}(\bar{k}) = \mt{l}(i)$, and add $k$ and $\bar{k}$ to $\mt{D}$
				\STATE \quad \quad \textbf{end for}
			\STATE \quad \textbf{end if} 
			\STATE \quad Remove $i$ from $\mt{D}$
		\ENDWHILE
		\STATE \textbf{for} each $i \in \widetilde{\mcl{I}}^j$ for each $j \in M$ \textbf{do}
			\STATE \quad add $i^{\mt{l}(i)}$ to $\mcl{I}_j$
		\STATE \textbf{end for}
		\STATE \textbf{for} each $i \in \widetilde{\mcl{I}}$ \textbf{do}
			\STATE \quad add $i^{\mt{l}(i)}$ to $\bar{\mcl{I}}$
		\STATE \textbf{end for}
		\STATE \textbf{for} each $i \in \widetilde{\mcl{J}}$ \textbf{do}
			\STATE \quad \textbf{if} $i^{\mt{l}(i)} = +$ \textbf{then}
				\STATE \quad \quad add $i$ to $\mcl{J}$
			\STATE \quad \textbf{else}	
				\STATE \quad \quad add $i$ to $\bar{\mcl{J}}$
			\STATE \quad \textbf{end if}			
		\STATE \textbf{end for}
	\end{algorithmic}
\end{algorithm}

%\begin{remark} \label{rem:primal-ECR-superset}
%	We note that the collection of all \ECR assignments obtained from Theorem~\ref{thm:primal-forest-converse} contains those corresponding to pairwise cancellation as a subset. 
%	More specifically, it is possible to obtain an \ECR assignment from Theorem~\ref{thm:primal-forest-converse} that may contain bilinear terms in the aggregated inequality that is not canceled during aggregation while appearing in two constraints.
%	See Example?? for an illustration of this case.
%\end{remark}	

\medskip	
In view of Theorem~\ref{thm:primal-forest-converse}, once we identify a forest structure with the desired conditions, we can use the steps in Algorithm~\ref{alg:primal-ECR} to determine the weight of each constraint in the corresponding \ECR assignment by following a path that starts from the arc associated with the base equality and reaches the node or arc associated with that constraint.
We illustrate this approach in the following example.

\medskip
\begin{example} \label{ex:primal-forest}
	Consider set $\mcl{S}$ with $m = 2$ and $\Xi$ that represents the primal network model corresponding to the graph $\mr{G} = (\mr{V},\mr{A})$ shown in Figure~\ref{fig:primal-tree}. 
	Similarly to Example~\ref{ex:primal-tree}, we refer to each arc in this network as a pair $(i,j)$ of its tail node $i$ and its head node $j$, and denote its corresponding flow variable as $x_{i,j}$.
	Assume that we are interested in finding \ECR assignments for class-$l^{+}$ where the base equality $l$ contains the bilinear term $y_1x_{1,5}$, i.e., $i' = (1,5)$ and $j' = 1$.
	According to Theorem~\ref{thm:primal-tree-converse}, we need to identify a forest structure that satisfies the conditions (i)--(iii) of Propositions~\ref{prop:primal-forest} and \ref{prop:primal-forest-pairwise}.
	In the parallel network $\mr{G}^1$, we select the forest $\bar{\mr{F}}^1$ composed of the tree $\bar{\mr{T}}^1_1$ with the node set $\{1, 2, 6\}$ and the tree $\bar{\mr{T}}^1_2$ with the node set $\{8\}$.
	In the parallel network $\mr{G}^2$, we select the forest $\bar{\mr{F}}^2$ composed of the tree $\bar{\mr{T}}^2_1$ with the node set $\{1, 4\}$.
	Therefore, we can form the set $\widetilde{\mcl{I}}^1 = \{1, 2, 6, 8\}$ and $\widetilde{\mcl{I}}^2 = \{1, 4\}$.
	We select the connection node set $\widetilde{\mcl{I}} = \{3\}$, and the connection arc set $\widetilde{\mcl{J}} = \{(8,4)\}$.
	It is easy to verify that these sets satisfy the conditions (i)--(iii) of Propositions~\ref{prop:primal-forest} and \ref{prop:primal-forest-pairwise}.
	Next, we determine the label of each node and arc in the above sets through applying Algorithm~\ref{alg:primal-ECR}.
	According to line 2 of this algorithm, we set $\mt{l}(1,5) = +$ in parallel network $\mr{G}^1$, and select $k = t(1,5) = 1 \in \widetilde{\mcl{I}}^1$.
	It follows from line 5 of the algorithm that $\mt{l}(1) = -$ and $k$ is added to $\mt{D}$.
	Following lines 10--13, we obtain for $\widetilde{\mcl{I}}^1$ that $\mt{l}(2) = \mt{l}(6) = -$, and for $\widetilde{\mcl{I}}$ that $\mt{l}(3) = +$.
	Then, from lines 26--28 of Algorithm~\ref{alg:primal-ECR}, we deduce for $\widetilde{\mcl{I}}^2$ that $\mt{l}(4) = -$, and from lines 11--13 for $\widetilde{\mcl{I}}^2$, we obtain that $\mt{l}(1) = -$.
	Lines 32--34 imply that $\mt{l}(8,4) = -$ for $\widetilde{\mcl{J}}$.
	Lastly, we conclude from lines 38--40 that $\mt{l}(8) = -$ for $\widetilde{\mcl{I}}^1$.
	As a result, following lines 47--59 of the algorithm, we obtain the \ECR assignment $\big[\{1^-, 2^-, 6^-, 8^-\},\{1^-, 4^-\},\{3^+\} \big| \emptyset, \{(8,4)\}\big]$ for class-$l^+$.
	Based on this assignment, we multiply the negative flow-balance constraints at nodes $1, 2, 6, 8$ with $y_1$, we multiply the negative flow-balance constraints at nodes $1, 4$ with $y_2$, we multiply the positive flow-balance constraint at node $3$ with $1-y_1-y_2$, and we multiply the upper bound constraint on variable $x_{8,4}$ with $1-y_1-y_2$, and aggregate them with the base bilinear equality corresponding to arc $(1,5)$ with weight 1 to obtain the aggregated inequality
	\begin{multline*}
		-z_{1,5} - y_1 x_{4,5} - y_1 x_{3,7} + y_1 x_{4,3} - y_2 x_{1,5} + y_2 x_{4,5} + y_2 x_{2,3} - y_2 x_{3,7}\\
		+ (f_1 + f_2 + f_3 + f_6 + f_8 - u_{8,4})y_1 + (f_1 + f_3 + f_4 - u_{8,4}) y_2 \\
		+ x_{3,7} -x_{2,3} - x_{4,3} - x_{8,4} - f_3 + u_{8,4} \geq 0,
	\end{multline*}
	where $f_i$ denotes the supply/demand value at node $i$, and $u_{i,j}$ denotes the upper bound for variable $x_{i,j}$.
	Following the relaxation step in the \ECR procedure, we may relax each of the seven remaining bilinear terms into two possible linear expressions, leading to 128 total \ECR inequalities.
	If implemented inside of a separation oracle, we can use Remark~\ref{rem:1-separation} to find the most violated inequality among these 128 inequalities efficiently in linear time.
\hfill	$\blacksquare$	
\end{example}
	
\medskip
We conclude this section with a remark on the practical implementation of the proposed \ECR inequalities.	
While there is an efficient separation algorithm to find a separating \ECR inequality among those created from a given \ECR assignment as noted in Remark~\ref{rem:1-separation}, the choice of the class of an \ECR assignment and its possible forest structure in the underlying network can lead to a large pool of candidates to consider during a branch-and-cut approach. 
Note that each \ECR inequality is obtained through an aggregation of the constraints of $\mcl{S}$ with proper weights.
In particular, given an \ECR assignment $\big[\mcl{I}_1,\dotsc,\mcl{I}_m,\bar{\mcl{I}} \big| \mcl{J},\bar{\mcl{J}}\big]$ for class-$l^{\pm}$, we aggregate the base inequality of the form $f_l(\vc{x},\vc{y},\vc{z}) \geq 0$ with constraints of the general form $h(\vc{y})g(\vc{x}) \geq 0$, where $h(\vc{y})$ represents the aggregation weight that could be $y_j$ or $1-\sum_{j \in M} y_j$, and where $g(\vc{x}) \geq 0$ denotes a linear network flow constraint that could be the flow-balance or variable bound constraints.
In most branch-and-cut approaches, the starting relaxation of the problem contains all linear side constraints on $\vc{x}$ and $\vc{y}$.
It follows that an optimal solution $(\bar{\vc{x}};\bar{\vc{y}};\bar{\vc{z}})$ of such relaxation that is to be separated satisfies $h(\bar{\vc{y}})g(\bar{\vc{x}}) \geq 0$ for all valid choices of function $h(\vc{y})$ and constraint $g(\vc{x}) \geq 0$.
Therefore, for the resulting aggregated inequality to be violated at a point $(\bar{\vc{x}};\bar{\vc{y}};\bar{\vc{z}})$, we must have the base inequality violated at that point, i.e, $f_l(\bar{\vc{x}},\bar{\vc{y}},\bar{\vc{z}}) < 0$.
This observation can be used to select the class and sign of the \ECR assignment to be generated during a separation process.
To this end, we may sort the \textit{residual} values defined by $\Psi_k = |\bar{y}_j\bar{x}_i - \bar{z}_k|$ for all $(i,j,k)$ such that $A^k_{ji} = 1$, and choose class $k$ as that associated with largest $\Psi_k$ with the class sign $+$ if $\bar{y}_j\bar{x}_i - \bar{z}_k < 0$, and class sign $-$ otherwise.
This perspective can shed light on explaining the observation that the \ECR inequalities obtained from \textit{fewer} aggregations tend to be more effective in practice as noted in \cite{davarnia:ri:ta:2017} and also observed in our experiments in Section~\ref{sec:computation}.
Specifically, the addition of constraints $h(\vc{y})g(\vc{x}) \geq 0$ in the aggregation can increase the left-hand-side value in the aggregated inequality when $h(\bar{\vc{y}})g(\bar{\vc{x}}) > 0$, which could reduce the odds of obtaining a violated aggregated inequality.

\medskip
Another observation that can be helpful for choosing the forest structures is considering the relaxation step in the \ECR procedure. 
As described therein, each remaining bilinear term $y_jx_i$ can be relaxed using either the bound constraints or the bilinear constraints.
The former case is equivalent to aggregating the inequality with a constraint of the form $h(\vc{y})g(\vc{x}) \geq 0$ where $h(\vc{y}) = y_j$ and $g(\vc{x}) \in \{x_i \geq 0, u_i-x_i \geq 0\}$, for which the previous argument holds about achieving a violation.
For the latter case, on the other hand, we aggregate the inequality with a bilinear constraint of the form $\pm(y_jx_i - z_k) \geq 0$ for $(i,j,k)$ such that $A^k_{j,i} = 1$, which can potentially lead to a violation depending on the value of $\Psi_k = |\bar{y}_j\bar{x}_i - \bar{z}_k|$.
As a result, we might choose forest structures that contain the nodes incident to arcs $i \in \mr{A}$ corresponding to the most violated values in $|\bar{y}_j\bar{x}_i - \bar{z}_k|$.
In our computational experiments presented in Section~\ref{sec:computation}, we use the above-mentioned approaches in our separation oracle to select the class of \ECR assignments and their forest structures, which show promising results.

%%%%%%%%%%%%%%%%%%%%%%%
%%%%%%%%%%%%%%%%%%%%%%%%%%%%%%%%%%%%%%%%%%%%%%
%%%%%%%%%%%%%%%%%%%%%%%%%%%%%%%%%%%%%%%%%%%%%%%%%%%%%%%%%%%%%%%%%%%%%
%%%%%%%%%%%%%%%%%%%%%%%%%%%%%%%%%%%%%%%%%%%%%%%%%%%%%%%%%%%%%%%%%%%%%%%%%%%%%%%%%%%%%%%%%%%%
\section{Computational Experiments} \label{sec:computation}
%%%%%%%%%%%%%%%%%%%%%%%%%%%%%%%%%%%%%%%%%%%%%%%%%%%%%%%%%%%%%%%%%%%%%%%%%%%%%%%%%%%%%%%%%%%%
%%%%%%%%%%%%%%%%%%%%%%%%%%%%%%%%%%%%%%%%%%%%%%%%%%%%%%%%%%%%%%%%%%%%%
%%%%%%%%%%%%%%%%%%%%%%%%%%%%%%%%%%%%%%%%%%%%%%
%%%%%%%%%%%%%%%%%%%%%%%

In this section, we present preliminary computational results to evaluate the impact of the cutting planes generated through the results of Section~\ref{sec:primal} on improving the relaxation gap for two network flow applications that contain bilinear constraints.
For these experiments, the codes are written in Python 3.8.8. and the optimization problems are solved using Gurobi 9.5.2 at its default settings.

%%%%%%%%%%%%%%%%%%%%%%%%%%%%%%%%%
%%%%%%%%%%%%%%%%%%%%%%%%%%%%%%%%%%%%%%%
\subsection{Fixed-Charge Network Flow Problem} \label{subsec:fixed-charge}
%%%%%%%%%%%%%%%%%%%%%%%%%%%%%%%%%%%%%%%
%%%%%%%%%%%%%%%%%%%%%%%%%%%%%%%%%
We study the bilinear formulation for the fixed-charge network flow problem proposed in \cite{rebennack:na:pa:2009}.
Consider a directed bipartite network structure where the nodes in the first partition are supply nodes with supply $s_i$ for $i \in S$, and the nodes in the second partition are demand nodes with demand $d_j$ for $j \in D$.
For each arc $(i,j)$, let $x_{ij}$ represent its flow and $u_{ij}$ denote its capacity.
We assume that the cost of sending one unit flow through arc $(i,j)$ is $0$ when $x_{ij} = 0$, and it is $c_{i,j} + t_{ij}x_{ij}$ when $0 < x_{ij} \leq u_{ij}$ for some constant $c_{ij} > 0$ and slope $t_{ij} > 0$.
In \cite{rebennack:na:pa:2009}, a piecewise concave underestimator of this cost function is considered, which calculates the cost as $c_{ij} + \frac{t_{ij}}{\epsilon_{ij}} x_{ij}$ when $0 \leq x_{ij} \leq \epsilon_{ij}$ for some $0 \leq \epsilon_{ij} \leq u_{ij}$, and $c_{i,j} + t_{ij}x_{ij}$ when $\epsilon_{ij} \leq x_{ij} \leq u_{ij}$.
Through introducing a binary variable $y_{ij} \in \{0,1\}$ that indicates the line segment in the piecewise function above, \cite{rebennack:na:pa:2009} proposes the following bilinear formulation for the corresponding fixed-charge problem.
\begin{align}
	\min \quad &\sum_{i \in S} \sum_{j \in D} \left( \left(c_{ij} + \frac{t_{ij}}{\epsilon_{ij}}\right)x_{ij} + t_{ij} y_{ij} - \frac{t_{ij}}{\epsilon_{ij}} z_{ij} \right) & \label{eq:fc-1}\\
	& \sum_{j \in D} x_{ij} \leq s_i, & \forall i \in S \label{eq:fc-2} \\
	& \sum_{i \in S} x_{ij} \geq d_j, & \forall j \in D  \label{eq:fc-3} \\
	& x_{ij} y_{ij} = z_{ij}, & \forall i \in S, j \in D  \label{eq:fc-4}\\
	& \sum_{i \in S} \sum_{j \in D} y_{ij} \leq b \label{eq:fc-5}\\
	& 0\leq x_{ij} \leq u_{ij},  & \forall i \in S, j \in D \label{eq:fc-6}\\
	& y_{ij} \in \{0,1\}, & \forall i \in S, j \in D \label{eq:fc-7}
\end{align}    
where the objective function computes the total cost of transportation.
Constraints \eqref{eq:fc-2} and \eqref{eq:fc-3} represent the flow-balance equations for the supply and demand nodes, respectively.
Constraint \eqref{eq:fc-4} defines the bilinear terms in the objective function.
Further, Constraint \eqref{eq:fc-5} imposes a budget requirement on the $\vc{y}$ variables.
%Note that the above formulation is defined over the continuous domain for $\vc{y}$ variables; see \cite{rebennack:na:pa:2009} for a detailed account.   
The base relaxation for this problem, which we refer to as the \textit{McCormick relaxation}, is obtained by replacing Constraints \eqref{eq:fc-4} with their McCormick relaxations and relaxing $\vc{y}$ variables to be continuous.
For each variable $y_{ij}$, the structure of the corresponding relaxation \eqref{eq:fc-2},\eqref{eq:fc-3},\eqref{eq:fc-4},\eqref{eq:fc-6} conforms to that of $\mcl{S}$ with $m = 1$.
Therefore, our goal is to assess the effectiveness of the \ECR inequalities associated with the tree structures in the above network using the results of Section~\ref{subsec:primal-single}.

\medskip
For our computational experiments, we generate 10 random instances for each size category based on the following data specifications.
We consider four different size categories by choosing the number of nodes in the complete bipartite network from $\{50, 100\}$ and by selecting the fraction value $\frac{\epsilon_{ij}}{u_{ij}}$ for arcs $(i,j)$ that have fixed-charge structure (which determines the position of the break point in the piecewise linear cost function) from $\{0.2, 05\}$. 
The number of $y$ variables is equal to 20\% of the total number of arcs, where each such variable is associated with a randomly selected arc. Supply and demand parameters are chosen from a discrete uniform distribution between $(20,50)$. The arc capacity is produced from a discrete uniform distribution between $(1,50)$, and the budget $b$ is set at 20\% of the total number of $y$ variables. The slope $t$ for each cost function is randomly generated from a discrete uniform distribution between $(50,100)$. If an arc $(i,j)$ follows the fixed-charge structure, the value of $c_{ij}$ is selected from a uniform distribution between $(1,5)$. Otherwise, this value is generated from a uniform distribution between $(10,20)$.

\medskip
Table~\ref{tab:fixed-charge} shows the results of adding the \ECR cutting planes for the single-variable relaxations of \eqref{eq:fc-2}--\eqref{eq:fc-7} as described above.
The first column contains the network size (the number of nodes), the second column shows the fraction $\frac{\epsilon_{ij}}{u_{ij}}$ of arcs $(i,j)$ that have fixed-charge structure (which determines the position of the break point in the piecewise linear cost function), and the third column indicates the instance number.
The fourth column shows the optimal value of the McCormick relaxation of \eqref{eq:fc-1}--\eqref{eq:fc-7}.
The fifth and sixth columns, respectively, represent the optimal value and the solution time (in seconds) of \eqref{eq:fc-1}--\eqref{eq:fc-7} obtained by the Gurobi solver at its default setting.
If an optimal solution is not found within the time limit of 1000 seconds, the best solution found at termination is reported, which is indicated by ``$\geq 1000$" in the ``Time" column.
The percentage of the gap improvement obtained by the Gurobi solver at the root node (compared to the McCormick relaxation), before the start of branch-and-bound, is given in the seventh column.
The next two columns under ``Full \ECR" contain the result of adding all violated \ECR inequalities obtained from tree structures according to Theorem~\ref{thm:primal-tree-converse} with up to two cancellations (i.e., three aggregations).
These cuts are added in loops after the LP relaxation is solved to separate the current optimal solution until the improvement in the optimal value is less than 1\%.
To find the most violated \ECR inequalities produced from an \ECR assignment, we use the technique in Remark~\ref{rem:1-separation}.
The column ``Gap" contains the gap improvement obtained by adding these \ECR inequalities compared to the optimal value of the McCormick relaxation reported in column 6.
The next column shows the total solution time to add these inequalities.
The column ``Gap" under ``Separation \ECR" includes the result of adding the above \ECR inequalities through a separation oracle according to that discussed in Section~\ref{sec:primal}.
In particular, for a current optimal solution $(\bar{\vc{x}}; \bar{\vc{y}}; \bar{\vc{z}})$, we consider the \ECR assignment class and sign associated with the 35 largest values for $\Psi_{ij} = |\bar{x}_{ij}\bar{y}_{ij} - \bar{z}_{ij}|$ for all $(i,j)$.
We add the resulting \ECR inequalities in loops as discussed above.
The next column shows the solution time when using this separation method.  
The last two columns under ``RLT" contain the results of using the reformulation-linearization technique (RLT) of level 1 \cite{sherali:ad:1990}.
Specifically, the column ``Full RLT" shows the gap closure achieved by augmenting the McCormick relaxation with the linearized constraints as a result of multiplying each constraint of the formulation with $y_{ij}$ and $(1-y_{ij})$ for all $i \in S$ and $j \in D$.
The next two columns show the gap improvement obtained by this method compared to the McCormick relaxation, and the solution time for this approach.
The symbol ``-" indicates that the model has not been solved to optimality within the time limit of 5000 seconds.
%The column ``Separation RLT" includes the dual bound obtained by solving cut-generating linear programs (CGLPs) that produce cutting planes through the projection of the RLT formulation onto the original space of variables.
%These cutting planes are added in loops in a manner similar to that described above.
%The last two columns contains the gap improvement and the solution time for this approach. 
The last row for each problem size reports the average values over the 10 random instances.

\medskip
The results in Table~\ref{tab:fixed-charge} indicate the effectiveness of the \ECR approach in five areas.
First, the gap improvement values show the effectiveness of the proposed \ECR inequalities based on the tree structures in improving the gap closure and strengthening the classical McCormick relaxation.
Second, the gap improvement obtained through addition of the \ECR inequalities exceeds that obtained by the Gurobi solver at the root node.
Third, these results also imply that the \ECR approach is more effective in improving the dual bounds compared to the RLT application; the \ECR cuts provide similar gap improvement values in a much smaller time. 
Fourth, they support the general observation for aggregation-based methods, such as the \ECR, that the inequalities obtained from fewer aggregations (up to three in these experiments) tend to be the most impactful, as they account for a considerable portion of the total gap closure for most instances.
This is evident from comparing the improvement achieved by the \ECR cuts with that of the RLT, which provides an upper bound for the possible gap improvement that can be obtained from single-variable disjunctive relaxations of the problem. 
Fifth, these experiments demonstrate the effectiveness of the proposed separation method, which achieves similar gap improvement values in much smaller time compared to the case without separation.
These observations show promise for an efficient implementation of the \ECR technique to solve practical problems.

\begin{table}[!ht]
	\centering
	\caption{Evaluating \ECR cutting planes for the fixed-charge network flow problem}
	\label{tab:fixed-charge}
	\resizebox{1.00\textwidth}{!}{
		%		\begin{minipage}{\textwidth}
			% Please add the following required packages to your document preamble:
			%\usepackage{multirow}
			\begin{tabular}{|c|c|c||c||c|c||c||c|c||c|c||c|c|}
				\hline
				Node \# & Frac. & \# & LP & \multicolumn{3}{c||}{Solver} & \multicolumn{4}{c||}{Tree \ECR} & \multicolumn{2}{c|}{RLT} \\ \cline{5-13} 
				& & & & Opt. & Time & Root & Full & Time & Sep. & Time & Full & Time \\ \hline
                    50 & 0.2 & 1 & 4437.82 & 5322.36 & 735.83 & 0.58 & 0.72 & 94.42 & 0.66 & 2.19 & 0.72 & 859.62\\
                     & & 2 & 3815.63 & 4617.38 & 128.63 & 0.69 & 0.82 & 90.74 & 0.8 & 2.58 & 0.83 & 1878.66\\
                     &  & 3 & 4769.89 & 5654.0 & 92.98 & 0.57 & 0.77 & 91.15 & 0.72 & 2.44 & 0.77 & 824.01\\
                     &  & 4 & 5361.71 & 5935.23 & 378.29 & 0.44 & 0.69 & 93.58 & 0.69 & 2.54 & 0.69 & 1913.67\\
                     &  & 5 & 4424.93 & 5446.63 & $\geq 1000$ & 0.65 & 0.76 & 104.65 & 0.73 & 2.66 & 0.77 & 2714.41\\
                     &  & 6 & 4508.74 & 5302.18 & 322.14 & 0.62 & 0.83 & 103.06 & 0.82 & 3.11 & 0.83 & 1186.32\\
                     &  & 7 & 4693.38 & 5590.65 & $\geq 1000$ & 0.56 & 0.73 & 72.07 & 0.7 & 1.97 & 0.74 & 2452.11\\
                     &  & 8 & 4441.68 & 5145.72 & 283.13 & 0.58 & 0.74 & 94.69 & 0.7 & 2.45 & 0.75 & 1009.55\\
                     &  & 9 & 4233.48 & 5233.15 & 22.76 & 0.84 & 0.93 & 90.82 & 0.88 & 1.81 & 0.93 & 1116.21\\
                     &  & 10 & 4897.29 & 5636.04 & 710.44 & 0.64 & 0.77 & 97.04 & 0.76 & 2.58 & 0.77 & 838.35\\
                    \hline
                    avg &  & & &  & 467.42 & 0.62 & 0.78 & 93.22 & 0.75 & 2.44 & 0.78 & 1479.29   \\
                    \hline
                    50 & 0.5 & 1 & 4735.93 & 5347.03 & 254.85 & 0.55 & 0.77 & 90.53 & 0.75 & 2.1 & 0.78 & 197.58\\
                     &  & 2 & 4274.4 & 4856.78 & 10.23 & 0.61 & 0.87 & 72.98 & 0.83 & 2.24 & 0.89 & 293.13\\
                     &  & 3 & 5026.89 & 5610.48 & 7.87 & 0.75 & 0.86 & 71.24 & 0.84 & 1.93 & 0.88 & 199.79\\
                     &  & 4 & 4630.15 & 5246.89 & 98.49 & 0.65 & 0.85 & 67.95 & 0.82 & 2.45 & 0.86 & 198.27\\
                     &  & 5 & 4036.38 & 4669.68 & 106.82 & 0.71 & 0.85 & 91.06 & 0.81 & 1.9 & 0.87 & 236.8\\
                     &  & 6 & 4474.54 & 5151.84 & 198.48 & 0.59 & 0.81 & 87.82 & 0.75 & 1.72 & 0.83 & 211.01\\
                     &  & 7 & 5176.22 & 5646.79 & 80.32 & 0.63 & 0.82 & 67.16 & 0.8 & 2.41 & 0.83 & 316.09\\
                     &  & 8 & 5084.2 & 5721.84 & $\geq 1000$ & 0.54 & 0.79 & 67.33 & 0.77 & 2.53 & 0.8 & 254.97\\
                     &  & 9 & 4523.52 & 5088.3 & 105.66 & 0.63 & 0.79 & 87.54 & 0.77 & 2.64 & 0.8 & 241.79\\
                     &  & 10 & 4735.93 & 5347.03 & 254.85 & 0.55 & 0.77 & 90.53 & 0.75 & 2.1 & 0.78 & 197.58\\
				\hline
                    avg & &  & &  & 186.7 & 0.63 & 0.83 & 79.6 & 0.8 & 2.26 & 0.84 & 238.89 \\
                    \hline
                    100 & 0.2 & 1 & 7088.91 & 9213.18 & $\geq 1000$ & 0.36 & 0.51 & 1162.19 & 0.46 & 8.4 & - & $\geq 5000$\\
                     &  & 2 & 7146.18 & 8960.31 & $\geq 1000$ & 0.48 & 0.63 & 1142.41 & 0.57 & 8.74 & - & $\geq 5000$\\
                     &  & 3 & 7248.75 & 9040.84 & $\geq 1000$ & 0.34 & 0.47 & 882.9 & 0.42 & 9.28 & - & $\geq 5000$\\
                     &  & 4 & 7692.17 & 9836.81 & $\geq 1000$ & 0.3 & 0.43 & 902.66 & 0.39 & 8.51 & - & $\geq 5000$\\
                     &  & 5 & 6802.61 & 8709.36 & $\geq 1000$ & 0.4 & 0.55 & 1144.94 & 0.48 & 8.57 & - & $\geq 5000$\\
                     &  & 6 & 6950.35 & 8862.86 & $\geq 1000$ & 0.47 & 0.61 & 1272.74 & 0.57 & 8.58 & - & $\geq 5000$\\
                     &  & 7 & 7052.04 & 8963.63 & $\geq 1000$ & 0.43 & 0.59 & 1151.14 & 0.52 & 8.37 & - & $\geq 5000$\\
                     &  & 8 & 7018.5 & 8894.89 & $\geq 1000$ & 0.49 & 0.62 & 1169.81 & 0.57 & 10.64 & - & $\geq 5000$\\
                     &  & 9 & 6732.16 & 8723.6 & $\geq 1000$ & 0.41 & 0.53 & 1132.3 & 0.48 & 8.49 & - & $\geq 5000$\\
                     &  & 10 & 7088.91 & 9213.18 & $\geq 1000$ & 0.36 & 0.51 & 1162.19 & 0.46 & 8.4 & - & $\geq 5000$\\
                    \hline
                    avg &  & & &  & $\geq 1000$ & 0.41 & 0.55 & 1111.73 & 0.49 & 8.87 &   & $\geq 5000$\\
                    \hline
                    100 & 0.5 & 1 & 7277.28 & 8572.76 & $\geq 1000$ & 0.6 & 0.78 & 1537.78 & 0.7 & 8.52 & - & $\geq 5000$\\
                     &  & 2 & 6991.21 & 8398.58 & $\geq 1000$ & 0.63 & 0.75 & 1388.98 & 0.68 & 8.1 & - & $\geq 5000$\\
                     &  & 3 & 7486.13 & 8604.12 & $\geq 1000$ & 0.52 & 0.72 & 1331.14 & 0.63 & 7.76 & - & $\geq 5000$\\
                     &  & 4 & 6768.0 & 8145.13 & $\geq 1000$ & 0.52 & 0.7 & 1370.79 & 0.61 & 7.93 & - & $\geq 5000$\\
                     &  & 5 & 7018.42 & 8413.21 & $\geq 1000$ & 0.58 & 0.73 & 1439.65 & 0.67 & 8.33 & - & $\geq 5000$\\
                     &  & 6 & 7139.52 & 8558.94 & $\geq 1000$ & 0.51 & 0.67 & 1470.46 & 0.6 & 7.89 & - & $\geq 5000$\\
                     &  & 7 & 7362.48 & 8822.47 & $\geq 1000$ & 0.51 & 0.65 & 1519.48 & 0.58 & 8.3 & - & $\geq 5000$\\
                     &  & 8 & 7232.95 & 8613.9 & $\geq 1000$ & 0.6 & 0.76 & 1493.97 & 0.72 & 10.38 & - & $\geq 5000$\\
                     &  & 9 & 6881.73 & 8112.89 & $\geq 1000$ & 0.59 & 0.77 & 1510.83 & 0.72 & 10.47 & - & $\geq 5000$\\
                     &  & 10 & 7277.28 & 8572.76 & $\geq 1000$ & 0.6 & 0.78 & 1537.78 & 0.7 & 8.52 & - & $\geq 5000$\\
                     \hline
                     avg & &  &   &  & $\geq 1000$ & 0.57 & 0.73 & 1445.3 & 0.66 & 8.89 & & $\geq 5000$ \\
                    \hline
			\end{tabular}
			%		\end{minipage}
	}
\end{table}

%%%%%%%%%%%%%%%%%%%%%%%%%%%%%%%%%
%%%%%%%%%%%%%%%%%%%%%%%%%%%%%%%%%%%%%%%
\subsection{Transportation Problem with Conflicts} \label{subsec:transportation}
%%%%%%%%%%%%%%%%%%%%%%%%%%%%%%%%%%%%%%%
%%%%%%%%%%%%%%%%%%%%%%%%%%%%%%%%%
In this section, we study a generalization of the Red-Blue Transportation Problem---which is a class of the Transportation Problem with Conflicts---presented in \cite{vancroonenburg:de:go:sp:2014}.
In particular, we consider a balanced set of supply and demand nodes that are arranged in a bipartite network structure, where there is a route (arc) between each supply node to each demand node.
Each supply node $i \in S$ has a supply $s_i$, and each demand node $j \in D$ has a demand $d_j$.
For each arc $(i,j)$, we refer to the flow sent through that arc as $x_{ij}$ and to the capacity of that arc as $u_{ij}$. 
The default cost of transportation using arc $(i,j)$ is $c_{ij}$.
In addition to the regular transportation means, there are multiple transportation services that could be used to facilitate the transfer of goods via different routes.
Binary variable $y_k$ represents whether or not service $k \in K$ is used.
The relative incentive/cost of using service $k \in K$ for each arc $(i,j)$ is denoted by $r_{ij}^k$.
In the Red-Blue Transportation Problem, due to certain logistic conflicts, some pairs of services could not be used simultaneously, i.e., $y_ky_l = 0$ for each pair of services $k$ and $l$ that are in the conflict set $C$.
This problem can be formulated as
\begin{align}
	\min \quad &\sum_{i \in S} \sum_{j \in D} \left(c_{ij}x_{ij} + \sum_{k \in K} r_{ij}^k z_{ij}^k \right) & \label{eq:tr-1}\\
	& \sum_{j \in D} x_{ij} \leq s_i, & \forall i \in S \label{eq:tr-2} \\
	& \sum_{i \in S} x_{ij} \geq d_j, & \forall j \in D  \label{eq:tr-3} \\
	& x_{ij} y_k = z_{ij}^k, & \forall i \in S, j \in D, k \in K  \label{eq:tr-4}\\
	& y_k + y_l \leq 1, & \forall (k,l) \in C  \label{eq:tr-5}\\
	& 0\leq x_{ij} \leq u_{ij},  & \forall i \in S, j \in D \label{eq:tr-6}\\
	& y_k \in \{0,1\}, & \forall k \in K \label{eq:tr-7}
\end{align}
where the objective function calculates the total cost of transportation.
Constraints \eqref{eq:tr-2} and \eqref{eq:tr-3} represent the flow-balance equations for the supply and demand nodes, respectively.
Constraint \eqref{eq:tr-4} defines the bilinear terms in the objective function that capture the incentive collected through using different transportation services.
The conflicts between the pair of transportation services in $C$ are imposed by Constraint \eqref{eq:tr-5}.
The base relaxation for this problem, which we refer to as the \textit{McCormick relaxation}, is obtained by replacing Constraints \eqref{eq:tr-4} with their McCormick relaxations and relaxing $\vc{y}$ variables to be continuous.

\medskip
For our experiments, we generate random generated instances with the following specifications. 
We consider four size categories by choosing the number of nodes in the complete bipartite network from $\{50, 100\}$, and selecting the number of transportation services from $\{20, 30\}$. 
Supply and demand parameters are generated from a discrete uniform distribution between $(100,200)$, and the arc capacities are generated from a discrete uniform distribution between $(1,25)$.
The number of conflicts is set at 10\% of all possible pairwise conflicts between the transportation services. The cost $c_{ij}$ for each arc $(i,j)$ is chosen randomly from a uniform distribution between $(20,40)$ multiplied by the number of transportation services for each problem set, and the relative incentive/cost $r_{ij}^k$ of using service $k$ is randomly chosen from a uniform distribution between $(-10,10)$.

\medskip
Table~\ref{tab:transportation} shows the results of adding the \ECR cutting planes to the McCormick relaxation described above.
The first column contains the network size (the number of nodes), the second column represents the number of available transportation services (i.e. the number of $y$ variables), and the third column indicates the instance number.
The values in the next four columns are defined similarly to those of Table~\ref{tab:fixed-charge}.
The time limit to solve the MIP reformulation of the problem by Gurobi is set to 5000 seconds. If the problem is not solved to optimality within the time limit, its best solution found at termination is reported, which is indicated by ``$\geq 5000$" in the ``Time" column.
The columns under ``\ECR" contain the gap improvement results obtained by adding different types of \ECR inequalities to the McCormick relaxation.
In particular, ``Tree Full" represents the gap improvement obtained by adding \ECR inequalities with up to three aggregations (two cancellations) generated from the one-variable relaxations of \eqref{eq:tr-2}--\eqref{eq:tr-7} where only one $y$ variable is considered.
For this approach, we use the \ECR results of Theorem~\ref{thm:primal-tree-converse} to identify the tree structures for each one-variable relaxation and add the resulting cutting planes for each relaxation separately through loops as previously described in Section~\ref{subsec:fixed-charge}.
The next column contains the time to implement these tree inequalities.
The column ``Forest \ECR" includes the gap improvement obtained by adding \ECR inequalities with up to three aggregations obtained for $\mcl{S}$ where two $y$ variables in a conflict are considered in their original simplex.
To add \ECR cutting planes, we consider the forest structures according to Theorem~\ref{thm:primal-forest-converse}.
The next column shows the time it takes to use this approach.
The next two columns ``Sep." and ``Time" indicate the gap closure and the solution time to implement these \ECR cutting planes using the separation oracle introduced in Section~\ref{subsec:primal-multi}, respectively.
The values in the columns with the header ``RLT" are defined similarly to those of Table~\ref{tab:fixed-charge}.

\medskip

It is evident from the results of Table~\ref{tab:transportation} that the \ECR cuts obtained from the forest structures that simultaneously consider multiple $y$ variables significantly outperform those obtained from the tree structures that consider $y$ variables individually, showing the effectiveness of the introduced class of \ECR inequalities with the pairwise cancellation property.
Further, these results show the remarkable impact of using a separation oracle to produce the \ECR inequalities on reducing the solution time, especially for larger size problems where the gap improvement by the \ECR inequalities based on forest structures outperform the RLT approach in both the gap improvement and solution time.

\begin{table}[!ht]
	\centering
	\caption{Evaluating \ECR cutting planes for the transportation problem with conflicts}
	\label{tab:transportation}
	\resizebox{1.00\textwidth}{!}{
		%		\begin{minipage}{\textwidth}
			% Please add the following required packages to your document preamble:
			%\usepackage{multirow}
			\begin{tabular}{|c|c|c||c||c|c||c||c|c||c|c|c|c||c|c|}
				\hline
				Node \# & $|K|$ & \# & LP & \multicolumn{3}{c||}{Solver} & \multicolumn{6}{c||}{\ECR} & \multicolumn{2}{c|}{RLT} \\ \cline{5-15} 
				& & & & Opt. & Time & Root & Tree Full & Tree Time & Forest Full & Forest Time & Sep. & Time & Full & Time \\ \hline
        50 & 20 & 1 & 1816461.35 & 1845306.7 & 58.71 & 0.14 & 0.24 & 106.05 & 0.54 & 424.88 & 0.5 & 42.55 & 0.63 & 41.39\\
         &  & 2 & 1911383.7 & 1932617.96 & 61.13 & 0.13 & 0.3 & 92.53 & 0.6 & 407.09 & 0.56 & 48.2 & 0.68 & 44.42\\
         &  & 3 & 1791541.33 & 1815825.43 & 66.51 & 0.03 & 0.27 & 107.84 & 0.56 & 449.59 & 0.53 & 50.07 & 0.62 & 37.26\\
         &  & 4 & 1916558.28 & 1938151.82 & 53.96 & 0.12 & 0.25 & 102.08 & 0.53 & 605.55 & 0.5 & 46.53 & 0.7 & 37.42\\
         &  & 5 & 1883973.93 & 1907566.87 & 102.68 & 0.13 & 0.25 & 105.91 & 0.54 & 443.17 & 0.52 & 43.61 & 0.65 & 36.89\\
         &  & 6 & 1754821.67 & 1779059.14 & 45.84 & 0.06 & 0.29 & 104.15 & 0.6 & 555.97 & 0.56 & 43.03 & 0.67 & 39.77\\
         &  & 7 & 1859891.92 & 1882935.97 & 61.43 & 0.07 & 0.33 & 90.97 & 0.62 & 609.47 & 0.58 & 44.93 & 0.7 & 37.78\\
         &  & 8 & 1878222.67 & 1907158.85 & 71.31 & 0.04 & 0.24 & 102.66 & 0.54 & 423.88 & 0.51 & 44.37 & 0.56 & 32.62\\
         &  & 9 & 1888464.09 & 1917359.49 & 104.42 & 0.13 & 0.27 & 110.78 & 0.54 & 633.7 & 0.51 & 49.51 & 0.63 & 34.44\\
         &  & 10 & 1828942.4 & 1853959.15 & 81.38 & 0.13 & 0.31 & 104.1 & 0.56 & 427.72 & 0.52 & 47.27 & 0.67 & 34.09\\
				\hline
        avg &  &  &  &  & 70.74 & 0.1 & 0.28 & 102.71 & 0.56 & 498.1 & 0.53 & 46.01 & 0.65 & 37.61\\
				\hline
        50 & 30 & 1 & 2735130.47 & 2769326.44 & 416.31 & 0.18 & 0.29 & 173.18 & 0.57 & 1382.33 & 0.56 & 87.34 & 0.69 & 107.66\\
         &  & 2 & 2758928.32 & 2797673.91 & 219.19 & 0.19 & 0.24 & 139.81 & 0.55 & 854.01 & 0.52 & 76.51 & 0.67 & 98.79\\
         &  & 3 & 2790718.79 & 2830074.49 & 152.08 & 0.16 & 0.22 & 141.19 & 0.56 & 1179.99 & 0.52 & 76.41 & 0.65 & 98.27\\
         &  & 4 & 2948401.96 & 2977573.06 & 75.43 & 0.22 & 0.37 & 148.31 & 0.69 & 1215.35 & 0.66 & 75.41 & 0.85 & 117.27\\
         &  & 5 & 2848728.26 & 2884380.56 & 189.93 & 0.16 & 0.27 & 153.52 & 0.58 & 1010.3 & 0.55 & 74.65 & 0.7 & 93.62\\
         &  & 6 & 2913904.34 & 2960520.37 & 510.74 & 0.24 & 0.2 & 177.89 & 0.53 & 1067.57 & 0.51 & 92.51 & 0.63 & 98.26\\
         &  & 7 & 2600592.36 & 2643997.19 & 367.13 & 0.16 & 0.26 & 184.94 & 0.57 & 1055.12 & 0.55 & 83.96 & 0.63 & 113.51\\
         &  & 8 & 2798975.64 & 2841262.04 & 324.84 & 0.2 & 0.25 & 176.37 & 0.53 & 1467.63 & 0.49 & 84.44 & 0.7 & 105.82\\
         &  & 9 & 2909192.25 & 2948513.71 & 219.4 & 0.2 & 0.28 & 187.26 & 0.59 & 950.87 & 0.56 & 93.77 & 0.72 & 106.4\\
         &  & 10 & 2611547.24 & 2656632.4 & 936.96 & 0.16 & 0.23 & 150.52 & 0.54 & 878.4 & 0.52 & 82.67 & 0.62 & 91.72\\
                    \hline
        avg &  &  &  &  & 341.2 & 0.19 & 0.26 & 163.3 & 0.57 & 1106.15 & 0.54 & 82.77 & 0.69 & 103.13\\
                    \hline
        100 & 20 & 1 & 3160875.04 & 3261525.49 & 2901.52 & 0.11 & 0.14 & 983.68 & 0.44 & 3932.78 & 0.42 & 409.97 & 0.41 & 1296.16\\
         &  & 2 & 3322767.42 & 3418349.48 & 1388.27 & 0.11 & 0.15 & 986.57 & 0.46 & 3821.98 & 0.44 & 328.39 & 0.44 & 1935.97\\
         &  & 3 & 3264069.46 & 3356160.8 & 1194.99 & 0.13 & 0.16 & 1032.29 & 0.49 & 3932.34 & 0.46 & 359.18 & 0.46 & 1344.05\\
         &  & 4 & 3301911.59 & 3398167.09 & 1182.85 & 0.11 & 0.13 & 910.19 & 0.47 & 3275.66 & 0.45 & 314.41 & 0.41 & 2110.84\\
         &  & 5 & 3334562.36 & 3435820.55 & 2770.33 & 0.11 & 0.14 & 974.61 & 0.44 & 3806.47 & 0.42 & 318.75 & 0.41 & 2052.78\\
         &  & 6 & 3316661.88 & 3414631.77 & 2088.83 & 0.09 & 0.16 & 960.33 & 0.46 & 3912.27 & 0.44 & 324.51 & 0.41 & 2661.02\\
         &  & 7 & 3341341.4 & 3438435.11 & 1519.65 & 0.03 & 0.16 & 1019.78 & 0.5 & 3956.82 & 0.48 & 347.23 & 0.37 & 1227.38\\
         &  & 8 & 3252602.14 & 3351139.62 & 2477.65 & 0.12 & 0.14 & 788.56 & 0.46 & 3174.85 & 0.43 & 313.43 & 0.42 & 1082.16\\
         &  & 9 & 3273121.61 & 3369045.93 & 1647.25 & 0.12 & 0.14 & 745.21 & 0.44 & 3025.82 & 0.42 & 287.02 & 0.43 & 1201.9\\
         &  & 10 & 3231657.29 & 3339983.21 & 2644.37 & 0.03 & 0.13 & 791.88 & 0.4 & 3242.4 & 0.38 & 333.17 & 0.33 & 993.24\\
                    
                    \hline
        avg &  &  &  &  & 1981.57 & 0.1 & 0.15 & 919.31 & 0.45 & 3608.14 & 0.44 & 333.6 & 0.41 & 1590.55\\
        \hline
        100 & 30 & 1 & 4834376.78 & 4968183.0 & $\geq 5000$ & 0.19 & 0.17 & 1866.58 & 0.55 & 10207.78 & 0.54 & 1179.07 & 0.51 & 2243.39\\
         &  & 2 & 4803153.48 & 4955801.77 & $\geq 5000$ & 0.19 & 0.12 & 1642.88 & 0.51 & 11944.64 & 0.49 & 994.29 & 0.46 & 2236.58\\
         &  & 3 & 4995871.11 & 5147567.66 & $\geq 5000$ & 0.13 & 0.14 & 1838.84 & 0.5 & 11362.07 & 0.5 & 1091.86 & 0.43 & 2322.86\\
         &  & 4 & 4961244.87 & 5143262.24 & $\geq 5000$ & 0.17 & 0.12 & 1976.34 & 0.46 & 12195.03 & 0.44 & 1308.95 & 0.42 & 2062.33\\
         &  & 5 & 4863008.91 & 5024866.15 & $\geq 5000$ & 0.17 & 0.13 & 1833.8 & 0.49 & 11701.45 & 0.48 & 984.07 & 0.45 & 2238.44\\
         &  & 6 & 4902278.58 & 5060625.13 & $\geq 5000$ & 0.17 & 0.13 & 1498.19 & 0.49 & 8708.22 & 0.49 & 1099.36 & 0.45 & 1983.21\\
         &  & 7 & 4866695.51 & 5017000.38 & $\geq 5000$ & 0.2 & 0.14 & 1517.02 & 0.52 & 8596.32 & 0.51 & 1331.52 & 0.49 & 1886.76\\
         &  & 8 & 5023334.69 & 5184894.21 & $\geq 5000$ & 0.22 & 0.13 & 1437.16 & 0.51 & 9807.23 & 0.5 & 1086.44 & 0.48 & 1701.46\\
         &  & 9 & 5043560.92 & 5198406.71 & $\geq 5000$ & 0.15 & 0.14 & 1420.03 & 0.49 & 8983.96 & 0.48 & 1382.81 & 0.45 & 1784.4\\
         &  & 10 & 4969879.76 & 5117213.58 & $\geq 5000$ & 0.18 & 0.14 & 1379.89 & 0.5 & 8421.21 & 0.49 & 1189.98 & 0.47 & 1911.2\\
        \hline
        avg &  &  &  &  & $\geq 5000$ & 0.18 & 0.14 & 1641.07 & 0.5 & 10192.79 & 0.49 & 1164.84 & 0.46 & 2037.06\\
        \hline
        
			\end{tabular}
			%		\end{minipage}
	}
\end{table}

%%%%%%%%%%%%%%%%%%%%%%%
%%%%%%%%%%%%%%%%%%%%%%%%%%%%%%%%%%%%%%%%%%%%%%
%%%%%%%%%%%%%%%%%%%%%%%%%%%%%%%%%%%%%%%%%%%%%%%%%%%%%%%%%%%%%%%%%%%%%
%%%%%%%%%%%%%%%%%%%%%%%%%%%%%%%%%%%%%%%%%%%%%%%%%%%%%%%%%%%%%%%%%%%%%%%%%%%%%%%%%%%%%%%%%%%%
\section{Conclusion} \label{sec:conclusion}
%%%%%%%%%%%%%%%%%%%%%%%%%%%%%%%%%%%%%%%%%%%%%%%%%%%%%%%%%%%%%%%%%%%%%%%%%%%%%%%%%%%%%%%%%%%%
%%%%%%%%%%%%%%%%%%%%%%%%%%%%%%%%%%%%%%%%%%%%%%%%%%%%%%%%%%%%%%%%%%%%%
%%%%%%%%%%%%%%%%%%%%%%%%%%%%%%%%%%%%%%%%%%%%%%
%%%%%%%%%%%%%%%%%%%%%%%

We study a bipartite bilinear set, where the variables in one partition belong to a network flow model, and the variables in the other partition belong to a simplex.
We design a convexification technique based on the aggregation of side constraints with appropriate weights, which produces an important class of facet-defining inequalities for the convex hull of the bilinear set, describing the convex hull for the special case where the simplex contains a single variable.
We show that each such inequality can be obtained by considering the constraints corresponding to the nodes of the underlying network that form a special tree or forest structure.
This property leads to an explicit derivation of strong inequalities through identifying special graphical structures in the network model. 
These inequalities can be added to the classical McCormick relaxation to strengthen the relaxation and improve the dual bounds, as corroborated in the preliminary computational experiments conducted on network problems in different application areas.

%%%%%%%%%%%%%%%%%%%%%%%
%%%%%%%%%%%%%%%%%%%%%%%%%%%%%%%%%%%%%%%%%%%%%%
%%%%%%%%%%%%%%%%%%%%%%%%%%%%%%%%%%%%%%%%%%%%%%%%%%%%%%%%%%%%%%%%%%%%%
%%%%%%%%%%%%%%%%%%%%%%%%%%%%%%%%%%%%%%%%%%%%%%%%%%%%%%%%%%%%%%%%%%%%%%%%%%%%%%%%%%%%%%%%%%%%
%\section*{Acknowledgments} 
%%%%%%%%%%%%%%%%%%%%%%%%%%%%%%%%%%%%%%%%%%%%%%%%%%%%%%%%%%%%%%%%%%%%%%%%%%%%%%%%%%%%%%%%%%%%
%%%%%%%%%%%%%%%%%%%%%%%%%%%%%%%%%%%%%%%%%%%%%%%%%%%%%%%%%%%%%%%%%%%%%
%%%%%%%%%%%%%%%%%%%%%%%%%%%%%%%%%%%%%%%%%%%%%%
%%%%%%%%%%%%%%%%%%%%%%%

%This work was supported in part by the AFOSR Grant FA9550-23-1-0183, and the NSF Grant CMMI-2338641.
%The authors thank the associate editor and the two anonymous referees for suggestions that helped improve the paper. 

%%%%%%%%%%%%%%%%%%%%%%%%%%%%%%
%%%%%%%%%%%%%%%%%%%%%%%%%%%%%%%%%%%%%%%%%%%%%%%%%%%%%%%%%%%%
%%%%%%%%%%%%%%%%%%%%%%%%%%%%%%%%%%%%%%%%%%%%%%%%%%%%%%%%%%%%%%%%%%%%%%%%%%%%%%%%%%%%%%%%%%
% References
%%%%%%%%%%%%%%%%%%%%%%%%%%%%%%%%%%%%%%%%%%%%%%%%%%%%%%%%%%%%%%%%%%%%%%%%%%%%%%%%%%%%%%%%%%
%%%%%%%%%%%%%%%%%%%%%%%%%%%%%%%%%%%%%%%%%%%%%%%%%%%%%%%%%%%%
%%%%%%%%%%%%%%%%%%%%%%%%%%%%%%

%\bibliographystyle{plain}
% \bibliographystyle{spbasic}
%\bibliography{../References/IPDD}
\bibliographystyle{informs2014} % outcomment this and next line in Case 1
\bibliography{multilinear,interdiction,mip,Mine}

\newpage
%%%%%%%%%%%%%%%%%%%%%%%
%%%%%%%%%%%%%%%%%%%%%%%%%%%%%%%%%%%%%%%%%%%%%%
%%%%%%%%%%%%%%%%%%%%%%%%%%%%%%%%%%%%%%%%%%%%%%%%%%%%%%%%%%%%%%%%%%%%%
%%%%%%%%%%%%%%%%%%%%%%%%%%%%%%%%%%%%%%%%%%%%%%%%%%%%%%%%%%%%%%%%%%%%%%%%%%%%%%%%%%%%%%%%%%%%
\section*{Appendix} 
%%%%%%%%%%%%%%%%%%%%%%%%%%%%%%%%%%%%%%%%%%%%%%%%%%%%%%%%%%%%%%%%%%%%%%%%%%%%%%%%%%%%%%%%%%%%
%%%%%%%%%%%%%%%%%%%%%%%%%%%%%%%%%%%%%%%%%%%%%%%%%%%%%%%%%%%%%%%%%%%%%
%%%%%%%%%%%%%%%%%%%%%%%%%%%%%%%%%%%%%%%%%%%%%%
%%%%%%%%%%%%%%%%%%%%%%%

In this section, we present a summary of the derivation of the \ECR procedure presented in Section~\ref{sec:ECR}.
We refer the reader to \cite{davarnia:ri:ta:2017} for the detailed exposition on these results.
Consider set $\mcl{S}$ defined in Section~1.
It can be shown that $\vc{y}$ components in the extreme points of $\mcl{S}$ are binary-valued; see Proposition~2.1 in \cite{davarnia:ri:ta:2017}.
As a result, the convex hull of $\mcl{S}$ can be obtained as a disjunctive union of a finite number of polytopes, each fixing $\vc{y}$ at an extreme point of the simplex $\Delta_m$.
A description for this convex hull can be obtained in a higher-dimensional space using the reformulation-linearization technique \cite{sherali:ad:dr:1998,sherali:al:1992} or disjunctive programming \cite{balas:1998} through addition of new variables, as shown below.

\begin{equation}
	\left.
	\begin{array}{lll}
		& \pm A^k_{j.} \vc{w}^j \mp v^j_k \geq 0,  &\forall (k, j) \in K \times M \\
		& \mp \left(z_k - \sum_{j \in M} v^j_k\right) \geq 0, &\forall k \in K \\
		& E\vc{w}^j \geq \vc{f} y_j, &\forall j \in M \\
		& E\left(\vc{x}-\sum_{j \in M}\vc{w}^j\right) \geq \vc{f}\left(1-\sum_{j \in M}y_j\right), \\
		& \vc{0} \leq \vc{w}^j \leq \vc{u} y_j, &\forall j \in M \\
		& \vc{0} \leq \vc{x}-\sum_{j \in M}\vc{w}^j \leq \vc{u}\left(1-\sum_{j \in M}y_j\right). 
	\end{array}
	\right. \label{eq_ef}
\end{equation}

\medskip
In the above description, variables $\vc{w}^j$ and $v^j_k$ can be viewed as $y_j \vc{x}$ and $y_j z_k$, respectively, and the equalities are formulated as pairs of inequalities of opposite directions.
Because the convex hull description in \eqref{eq_ef} contains additional variables, we can use polyhedral projection to obtain a convex hull description in the space of original variables. 

\medskip
Define the dual variables associated with the constraints of \eqref{eq_ef} by $\vc{\alpha}^{j\pm}, \vc{\beta}^\pm \in \Re^{\kappa}_+$, $\vc{\gamma}^{j}, \vc{\theta} \in \Re^{\tau}_+$, and $\vc{\eta}^j, \vc{\rho}^j, \vc{\lambda}, \vc{\mu} \in \Re^{n}_+$, respectively.
It follows from Proposition 2.3 of \cite{davarnia:ri:ta:2017} that the collection of inequalities
\begin{equation}
	\sum_{i \in N} q_i(\vc{\pi}) x_i + \sum_{j \in M} r_j(\vc{\pi}) y_j + \sum_{k \in K} s_k(\vc{\pi}) z_k \ge t(\vc{\pi}), \label{eq:proj-facet}
\end{equation}
where
\begin{equation*}
	\begin{array}{rll}
		q_i(\vc{\pi}) &=& \sum_{t \in T}E_{ti} \theta_t + \lambda_i - \mu_i \nonumber \\
		r_j(\vc{\pi}) &=& \sum_{t \in T}f_t \left(\theta_t - \gamma^j_t\right) + \sum_{i \in N}\left(\rho^j_i - \mu_i\right) \nonumber \\
		s_k(\vc{\pi}) &=& -\left(\beta^+_k - \beta^-_k\right) \nonumber\\
		t(\vc{\pi}) &=& \sum_{t \in T}f_t \theta_t - \sum_{i \in N}\mu_i, \nonumber
	\end{array}
\end{equation*}

\medskip
for all extreme rays $\vc{\pi} = \left(\vc{\beta}^+;\vc{\beta}^-;\{\vc{\gamma}^j\}_{j \in M};\vc{\theta};\{\vc{\eta}^j\}_{j \in M};\{\vc{\rho}^j\}_{j \in M};\vc{\lambda};\vc{\mu}\right)$ of the projection cone
\begin{multline*}
	\mathcal{C} =
	\left\{
	\vc{\pi} \in \Re_+^{2\kappa + (m+1)(\tau + 2n)} \, \middle| \, \sum_{k \in K}A^k_{ji}\left(\beta^+_k - \beta^-_k\right) + \sum_{t \in T} E_{ti} \left(\gamma^j_t - \theta_t\right) \right. \\
	\left. + \eta^j_i - \rho^j_i - \lambda_i + \mu_i = 0, \, \, \,  \forall (i,j) \in N \times M \right\}
\end{multline*}
contains all \textit{non-trivial} facet-defining inequalities in the convex hull description of $\mcl{S}$.
A non-trivial inequality is one that cannot be implied by the linear constraints in the description of $\Xi$ and $\Delta_m$.
It is implied from Proposition 2.4 of \cite{davarnia:ri:ta:2017} that an extreme ray of $\mcl{C}$ that has components $\beta^{\pm}_k = 0$ for all $k \in K$ leads to a trivial inequality. 
Therefore, we may assume that $\beta_l^{\pm} = 1$ for some $l \in K$ in an extreme ray associated with a non-trivial facet-defining inequality of $\conv(\mcl{S})$ through proper scaling.
As a result, the search for the extreme rays of $\mcl{C}$ reduces to that of the extreme points of the restriction set $\mcl{C}^l$ for all $l \in K$, where
\begin{multline*} 
	\mathcal{C}^{l} =
	\left\{
	\vc{\pi}^l \in \Re^{2(\kappa-1) + (m+1)(\tau + 2n)}_{+} \, \middle|
	\sum_{k \in K_l}A^k_{ji}\left(\beta^+_k - \beta^-_k\right) + \sum_{t \in T} E_{ti} \left(\gamma^j_t - \theta_t\right) \right. \\
	\left. + \eta^j_i - \rho^j_i - \lambda_i + \mu_i = \pm A^l_{ji}, \, \forall (i,j) \in N \times M \right\},
\end{multline*}
where $K_l = K \setminus \{l\}$, and $\vc{\pi}^l$ is defined similarly to $\vc{\pi}$, but without elements $\beta^+_l$ and $\beta^-_l$.

\medskip
The components in the dual vector $\vc{\pi}^l$ can be interpreted as the weights used in the \ECR procedure as follows.
The fixing of the component $\beta^{\pm}_l$ is achieved in Step 1 of the \ECR procedure by picking a base equality $l$ with either $+1$ or $-1$ weights.
The components $\vc{\gamma}^j$ (resp. $\vc{\theta}$) can be viewed as the weights for $y_j$ (resp. $1 - \sum_{j \in M} y_j$) when multiplied with the non-bound constraints in $\Xi$ as demonstrated in Step 2 of the \ECR procedure.
Similarly, $\vc{\lambda}$ and $\vc{\mu}$ denote the weights for $1 - \sum_{j \in M} y_j$ when multiplied with the bound constraints in $\Xi$ in Step 3 of the \ECR procedure.
Finally, the relaxation step in the \ECR procedure will use the components $\vc{\beta}^{\pm}$ as the weights of the bilinear constraints in $K_l$, as well as the components $\vc{\eta}^j$ and $\vc{\rho}^j$ as the weights for $y_j$ when multiplied with the bound constraints in $\Xi$ to \textit{cancel} the remaining bilinear terms.
It can be shown that criteria (C1) and (C2) of the \ECR procedure provide necessary conditions for the selected weight vector to be an extreme point of $\mcl{C}^l$; see Theorem 2.7 in \cite{davarnia:ri:ta:2017}.
The resulting \ECR inequality is of the form \eqref{eq:proj-facet}, and the collection of all such inequalities contain all non-trivial facet-defining inequalities in $\conv(\mcl{S})$.

\medskip
In conclusion, we give an analogy that could better explain the difference between the aggregation procedure developed for the general set $\mcl{S}$ in \cite{davarnia:ri:ta:2017} as summarized here, and the explicit convexification results obtained for a special case of $\mcl{S}$ studied in this paper. Consider the disjunctive programming technique that provides a general procedure to obtain a convex hull formulation for the union of polytopes in a higher dimension. This procedure does not provide an explicit form for the valid inequalities of the convex hull in the original space of variables for a general case. When applied to a special case, where the polytopes have special structures (for example, obtained through split disjunction), one can design an algorithm to obtain the explicit form of the valid inequalities in the original space of variables (for example, split inequalities) that are more practical to generate and use. Here, one can view the aggregation procedure proposed in \cite{davarnia:ri:ta:2017} as the ``general disjuctive programming technique", and the explicit convexification results obtained for the network polytopes as the ``explicit results obtained for the special-structured polytopes".

\begin{comment}
	The components in the dual vector $\vc{\pi}^l$ can be interpreted as the weights used in the \ECR procedure as follows.
	The fixing of the component $\beta^{\pm}_l$ is achieved in Step 1 of the \ECR procedure by picking a base equality $l$ with either $+1$ or $-1$ weights.
	The components $\vc{\beta}^{\pm}$ represent the weights of the bilinear constraints in $K_l$ as described in Step 2 of the \ECR procedure.
	The components $\vc{\gamma}^j$ (resp. $\vc{\theta}$) can be viewed as the weights for $y_j$ (resp. $1 - \sum_{j \in M} y_j$) when multiplied with the non-bound constraints in $\Xi$ as demonstrated in Step 3 of the \ECR procedure.
	Similarly, $\vc{\lambda}$ and $\vc{\mu}$ denote the weights for $1 - \sum_{j \in M} y_j$ when multiplied with the bound constraints in $\Xi$ in Step 4 of the \ECR procedure.
	Finally, the relaxation step in the \ECR procedure will use the components $\vc{\eta}^j$ and $\vc{\rho}^j$ as the weights for $y_j$ when multiplied with the bound constraints in $\Xi$ to \textit{cancel} the remaining bilinear terms.
	It can be shown that criteria (C1) and (C2) of the \ECR procedure provide necessary conditions for the selected weight vector to be an extreme point of $\mcl{C}^l$; see Theorem 2.7 in \cite{davarnia:ri:ta:2017}.
	The resulting \ECR inequality is of the form \eqref{eq:proj-facet}, and the collection of all such inequalities contain all non-trivial facet-defining inequalities in $\conv(\mcl{S})$.
\end{comment}

\end{document}